\documentclass[12pt,leqno]{amsart}
\usepackage{geometry}
\usepackage{verbatim,upref,amsxtra,amsthm,xcolor,amssymb,amsxtra}
\usepackage[colorlinks,linkcolor=black,anchorcolor=blue,citecolor=blue]{hyperref}

\usepackage[mathscr]{eucal}
\usepackage{bbm}
\usepackage{dsfont}

\usepackage[initials,nobysame]{amsrefs}

\usepackage{mathtools}
\usepackage[normalem]{ulem}

\usepackage{varioref}

\newcommand{\Z}{\mathbb Z}
\newcommand{\N}{\mathbb N}
\renewcommand{\Re}{\mathsf{Re\, }}
\newcommand{\R}{\mathbb R}
\newcommand{\T}{\mathbb T}
\newcommand{\eps}{\varepsilon}
\newcommand{\bm}{{\boldsymbol m}}
\newcommand{\bn}{{\boldsymbol n}}
\newcommand{\bz}{{\boldsymbol z}}
\newcommand{\bk}{{\boldsymbol k}}
\newcommand{\bl}{{\boldsymbol \ell}}
\newcommand{\bw}{{\boldsymbol w}}
\newcommand{\bo}{{\boldsymbol 0}}

\renewcommand{\S}{{\mathcal S}}

\newcommand{\Dim}{\rm Dim}

\def\mybf #1{\textbf{\textit{#1}}}

\theoremstyle{plain}
\newtheorem{thm}{Theorem}[section]
\newtheorem{lem}[thm]{Lemma}
\newtheorem{prop}[thm]{Proposition}
\newtheorem{cor}[thm]{Corollary}

\theoremstyle{definition}
\newtheorem{defn}[thm]{Definition}

\newtheorem{exmp}[thm]{Example}
\theoremstyle{remark}

\newtheorem{Question}{\bf  Question}

\numberwithin{equation}{section}

\date{} 

\begin{document}


\frenchspacing

\setlength\marginparsep{8mm}
\setlength\marginparwidth{20mm}

\title[Khintchin conjecture and  related topics]{Khintchin conjecture and  related topics} 

\author{Aihua FAN}
\address{(A. Fan) 
LAMFA, UMR 7352 CNRS, University of Picardie, 33 rue Saint Leu, 80039 Amiens, France
and Wuhan Institute for Math \& AI, Wuhan University, Wuhan
430072, China}
\email{ai-hua.fan@u-picardie.fr}

\author{Shilei FAN}
\address{(S. Fan) School of Mathematics and Statistics, and Key Lab NAA--MOE, Central China Normal University, Wuhan 430079, China}
\email{slfan@ccnu.edu.cn}

\author{Herv\'e Queff\'elec}
\address{(H. Queff\'elec) CNRS, Laboratoire Paul Painlev\'e, UMR 8524 \& Labex CEMPI (ANR-LABX-0007-01),
Universit\'e de Lille                        Cit\'e Scientifique, B\^t. M2                           59655 Villeneuve d'Ascq Cedex,
FRANCE }
\email{herve.queffelec@univ-lille.fr}

\author{Martine Queff\'elec}
\address{(M. Queff\'elec) CNRS, Laboratoire Paul Painlev\'e, UMR 8524, \& Labex CEMPI (ANR-LABX-0007-01),
Universit\'e de Lille                          Cit\'e Scientifique, B\^t. M2                                  59655 Villeneuve d'Ascq Cedex,
FRANCE}
\email{martine.queffelec@univ-lille.fr}

\thanks{ The authors were supported by NSFC (grants No. 12331004, No.  12231013 and No. 12571205)}

	\begin{abstract} 
    Motivated by Khintchin’s 1923 conjecture, refuted by Marstrand in 1970, we  study the Khintchin class of functions associated to a given increasing sequence of integers.  When the Khintchin class contains $L^p(\mathbb{T})$, we call the sequence a $L^p$-Khintchin sequence. 
    We establish basic properties of Khintchin sequences, provide several constructions, and propose open problems for further research. We also initiate the study of Khintchin sequences of group endomorphisms on compact abelian groups. Under a Fourier-tightness assumption, we show that ergodicity (respectively, weakly mixing or strongly mixing) of a skew product of endomorphisms is equivalent to the corresponding property of the base system, supporting the idea that typical fiber orbits in such skew products should form Khintchin sequences. 
	
	\end{abstract}

    \maketitle
	
	\tableofcontents
	
	\section{Introduction}
	It is well known that for  any  sequence of distinct  integers $E=(\lambda_n)_{n\ge 1}\subset \Z$, 
	 the sequence of real numbers $(\lambda_n x)_{n\ge 1}$ is uniformly distributed modulo 1 (or $\!\!\!\mod 1$ for short) for Lebesgue almost every point $x \in [0,1)$.
 Consequently, for every {Riemann-integrable} function $f$ on the interval $[0,1)$,  identified with the circle $\mathbb{T}=\mathbb{R}/\mathbb{Z}$, we have

\begin{equation} \label{Riemann0} \lim_{N\to \infty}\  \frac{1}{N}\sum_{n=1}^N f(\lambda_{n} x)= \int_{\T} f d\mathbf{m}  \quad a.e.
\end{equation}
In 1923, Khintchin \cite{Khintchin1923} conjectured that when $E=\N$,  the property \eqref{Riemann0} would hold for bounded Borel measurable functions $f$. But
Marstrand \cite{Marstrand1970} refuted this conjecture in 1970 by a counterexample, which is the indicator function of some open set, and  Bourgain \cite{Bourgain1988Entropy} gave a non-constructive proof of Marstrand Theorem by his entropy method.
This leads us  to define the {{\it Khintchin class}} of  $E=(\lambda_n) \subset \N$ as 
\begin{equation}\label{khikhi} \mathcal{K}_E:=\Big\{f\in L^{1}(\T) :  \lim_{N\to \infty }\frac{1}{N}\sum_{n=1}^N f(\lambda_{n} x)= \int_{\T} f d\mathbf{m}  \hbox{\quad a.e.} \Big\}.
\end{equation}
The term {\em  Khintchin class}  was introduced by Fan-Queff\'elec-Queff\'elec  in \cite{AHM}, where the Furstenberg sequence (the increasing ordered sequence of the semigroup $\{2^n 3^m: n\ge 0, m\ge 3\}$) was investigated from different points of view.
If $\mathcal{K}_E$ contains some subspace $A$ of $L^{1}(\mathbb{T})$, we  say that $E$ is an {\em  $A$-Khintchin sequence} (it is termed $A^*$-sequence in \cite{Marstrand1970}).  
So, we will talk about {\em $L^p$-Khintchin sequences} for $1\le p\le \infty$.  An $L^1$-Khintchin sequence will be simply called {\em Khintchin sequence}. 
By the {\em  Khintchin problem}  we mean the determination of the Khintchin class of a given sequence or related problems. 
We state the typical question as follows: is a given sequence of integers a $L^p$-Khintchin sequence for some $1\le p\le \infty$\,?
The Khintchin problem is broad: there are many different integer sequences, and—as we will see—it also arises on any compact abelian group 
(See Section ~\ref{sect:groups}).
\medskip

Let us recall some known results. 
Raikov \cite{Raikov1936} proved that $(2^n)$ is a Khintchin sequence. This also follows from Birkhoff’s ergodic theorem, as pointed out by Riesz \cite{Riesz1945}.
 For the Furstenberg sequence $\{s_n\}$ which is the increasingly ordered sequence of the Furstenberg set $S:=\{2^n3^m: n\ge 0, m\ge 0\}$, based on Bewley's ergodic theorem of $\mathbb{N}^2$-actions Nair \cite{Nair1990} proved that
$(s_n)$ is a $L^1$-Khintchin sequence, after Marstrand had proved  that $(s_n)$ is a $L^\infty$-Khintchin sequence \cite{Marstrand1970}.  
There is an improvement of Marstrand's result about $S$ in \cite{AHM} and there is another proof of Nair's result in \cite{FFL}. 
However, Rosenblatt \cite{Rosenblatt1989} proved  that  the very lacunary sequence 
$(2^{2^n})$ is not  a $L^\infty$-Khintchin sequence.  Bourgain \cite{Bourgain1989}  proved that $(2^{n^2})$
 is a $L^p$-Khintchin sequence for every $p>1$, but Buczolich and Mauldin \cite{BM2010} proved that $(2^{n^2})$
 is not a $L^1$-Khintchin sequence. 
These results indicate that the Khintchin problem seems to be  subtle.  See the survey papers \cite{Buczolich2023} and \cite{C-Nair2014} for interesting related information.  
Let us point out that Beck \cite{Beck2015} ``saves'' Khinchin's conjecture in the continuous case by switching from a sequence to a continuous torus line. 
\medskip

In this paper, we study basic properties of Khintchin sequences and construct new examples via a multiplicative procedure. We also consider randomized versions of this procedure, leading to skew products whose ergodic, weakly mixing, and strongly mixing properties we analyze.

\medskip

In Section \ref{sect:integers}, we will  study basic properties of Khintchin sequences of integers, propose a procedure to produce Khintchin sequences and suggest some questions to be studied.  

Section ~\ref{sect:groups} is devoted to  Khintchin sequences of endomorphisms on compact abelian groups.

In Section ~\ref{sect:skewproduct}, we  study random sequences produced by skew products and  investigate the ergodic and  mixing properties of the resulting measure-preserving dynamical systems.

 \section{Khintchin sequences of integers}
 \label{sect:integers}

 In this section, we consider  sequences of integers which we view as sequences of endomorphisms of the compact group $\mathbb{T}=\mathbb{R}/\mathbb{Z}$. An investigation on an arbitrary compact abelian group will be undertaken in the next section. A great part of these  discussions on the special group $\mathbb{T}$ is also valid on general compact abelian groups. We are mainly interested in the almost everywhere convergence with respect to the Lebesgue measure $\mathbf{m}$ of the associated Khintchin operators acting on integrable functions. We begin by recalling the Banach principle and the Bellow-Jones principle, which serve as basic tools. Then the union, intersection, subsequence of Khintchin sequences are studied. Next we introduce a procedure  to produce possible  Khintchin sequences and the special multiplicative Thue-Morse sequence  is proved to be a Khintchin sequence. Many open questions are asked.

	
\subsection{Characterization of $L^p$-Khintchin sequences for $1\le p\le \infty$}
Our main concern is about the almost everywhere convergence of $T_nf(x)$ for every $f$ belonging to $L^p(\mathbb{T})$ (for some $1\le p\le \infty)$, where $T_n$'s are the operators, called {\em Khintchin operators} associated to $(\lambda_n)$,
$$
T_nf(x) =\frac{1}{n}\sum_{k=1}^n f(\lambda_k x), \quad  n\geq 1.
$$
Recall that $(\lambda_k)$ is a given sequence of integers. Our operators $T_n$'s are well defined for all Borel measurable functions $f$. They are all linear and contractive on each space $L^p(\mathbb{T})$ ($1\le p\le \infty$), because $x \mapsto \lambda_k x \mod 1$ preserves the Lebesgue measure $\mathbf{m}$ (we assume that $\lambda_k\not=0$). To answer the question of almost everywhere convergence for such sequence of operators, we shall apply  the Banach principle when $1\le p<\infty$, or apply  a modified Banach principle due to Bellow and Jones when $p=\infty$. 
 Let us first present  these two principles  in a general setting.
 Let us also point out that Sawyer's maximal inequality of weak type applies to our sequences of Khintchin operators (see \cite{Sawyer1966}). 
\medskip

Let $(X, \mathcal{B}, \mu)$ be a probability space. We denote by $L^0(X, \mathcal{B}, \mu)$ the space of all measurable functions.  As usual, $L^p(X, \mathcal{B}, \mu)$ ($1\le p\le \infty)$ denote the  Lebesgue spaces.  The space $L^0(X, \mathcal{B}, \mu)$ is a complete metric space with the metric
$$
d(f,g) = \int_X\frac{|f(x)-g(x)|}{1+|f(x)-g(x)|}d\mu(x).
$$
The convergence under this metric is the convergence in probability. 
\medskip

Let $(\mathbb{B}, \|\cdot\|)$ be a Banach space. We shall consider operators $S: \mathbb{B}\to L^0(X, \mathcal{B}, \mu)$ which are usually assumed linear or {\em sublinear}
i.e., $$|S (a f)| = |a|\, |S(f)|, \quad  |S(f+g)|\le |S(f)|+|S(g)|,
$$
and continuous, i.e. continuous in probability(namely, $L^0(X, \mathcal{B}, \mu)$ is equipped with the above metric). 
Let $(S_n)$ be a sequence of linear operators from $\mathbb{B}$ to $L^0(X, \mathcal{B}, \mu)$. For any $f\in \mathbb{B}$, we define the {\em maximal function}
$$
S^*f(x) = \sup_{n\ge 1} |S_n f(x)|.
$$
This maximal operator $S^*$ is sublinear.
In this setting, the classical Banach principle states as follows. 
\medskip

{\sc Banach Principle \cite{Banach1926}.} {\em Assume that each linear operator \[ S_n: (\mathbb{B}, \|\cdot\|) \to (L^0(X, \mathcal{B}, \mu), d(\cdot, \cdot))\] is  continuous. Suppose that $S^*f(x) <\infty$ a.e. for every $f\in \mathbb{B}$. Then
\begin{itemize}
\item[{\rm (B1)}] The maximal operator $S^*: \mathbb{B}\to L^0(X, \mathcal{B}, \mu)$ is  continuous at $0$;
\item[{\rm (B2)}] The set $\{f\in \mathbb{B}: \lim S_nf(x) \ {\rm exists} \ a.e. \}$ is closed in $\mathbb{B}$.
\end{itemize}
Conversely, if\ $\lim S_nf(x)$ exists a.e. for all $f\in \mathbb{B}$, it is trivially true that $S^*f(x) <\infty$ a.e. for every $f\in \mathbb{B}$. 
}
\medskip

 Remark that the continuity of the linear operators $S_n$, together with the assumption that
$S^{*}f(x) < \infty$ a.e.\ for every $f \in \mathcal{B}$, implies the existence of a nonincreasing
function $C : \mathbb{R}^{+} \to \mathbb{R}^{+}$ such that
\[
\lim_{t \to +\infty} C(t) = 0
\]
and
\[
\forall f \in \mathcal{B}, \ \forall \alpha > 0, \qquad
\mu\{x : S^{*}f(x) > \alpha \|f\|\} \le C(\alpha).
\]
See \cite[Theorem~5.1.2]{Weber09} for more details. The last inequality is usually referred to as a maximal inequality.

We have the following characterization for $L^p$-Khintchin sequences
($1\le p<\infty$).

\begin{thm}[Characterization of $L^p$-Khintchin sequences for $1\le p<\infty$] \label{thm:NSCp}
    Assume $1\le p<\infty$. A sequence $(\lambda_n)$ of integers is a $L^p$-Khintchin sequence if and only if
    $T^*f(x)<\infty$ $\mathbf{m}$-a.e. for every $f\in L^p(\mathbb{T})$, where $T^*f(x) = \sup_n |T_nf(x)|$ ($T_n$ being the Khintchin operator). 
\end{thm}
\begin{proof}
The necessity is trivial. 

Let us prove the sufficiency.
Apply the Banach principle to the Khintchin operators $(T_n)$ defined above with $\mathbb{B}=L^p(\mathbb{T}, \mathcal{B}(\mathbb{T}), \mathbf{m})$, $\mathbf{m}$ being the Lebesgue measure.  Since $1\le p<\infty$, the continuous functions are dense in $L^p(\mathbb{T})$. On the other hand, the limit $\lim_{n\to \infty}T_nf(x)$ exists a.e. for every continuous function $f$, because $(\lambda_k x)$ is  equidistributed  for almost every $x$.  Then, by the Banach principle, the limit $\lim_{n\to \infty }T_nf(x)$ exists a.e. for every $f\in L^p(\mathbb{T})$. To confirm that the limit is equal to $\int_\mathbb{T}f(x) dx$, we apply Theorem \ref{thm:L^p-conv} which ensures the $L^p$-convergence of $T_nf$ to $\int_\mathbb{T}f(x) dx$.
Theorem \ref{thm:L^p-conv} will be proved in the setting of group endomorphisms on a compact abelian group.
\end{proof}

This characterization of $L^p$-Khintchin sequence for $1\le p<\infty$ will be very useful.
But the above Banach principle is not effective when $p=\infty$ for  Khintchin problem. Because the continuous functions are not dense in $L^\infty(\mathbb{T})$. The Riemann integrable functions are not dense either. 
\medskip

Now let us present Bellow-Jones principle.  
Let $$
B_\infty=\{f\in L^\infty(\mu): \|f\|_\infty\le 1\}$$ be the closed unit ball  of $L^\infty(X, \mathcal{B}, \mu)$. It is known that the metric $d(f,g)$
restricted on $B_\infty$ is equivalent to the metric $\|f-g\|_p$ for every $1\le p<\infty$. So, $B_\infty$ is a complete metric space for any such metric defined by $d(\cdot, \cdot)$ or by $\|\cdot\|_p$. In the following, we always use this topology of $B_\infty$.
Bellow and Jones modified the Banach principle as follows.   
\medskip

{\sc Bellow-Jones Principle \cite{BJ1996}.}
{\em Assume that each linear operator 
 \[S_n: L^\infty(X, \mathcal{B}, \mu) \to L^0(X, \mathcal{B}, \mu)\] is  continuous,  i.e. continuous in probability. Suppose that  $S^*f(x) <\infty$ a.e. for every $f\in L^\infty(X, \mathcal{B}, \mu)$ and $S^*: B_\infty \to L^0(X, \mathcal{B},\mu)$ is continuous at $0$. Then
\begin{itemize}
\item[{\rm (B3)}] The set $\{f\in B_\infty: \lim S_nf(x) \ {\rm exists} \ a.e. \}$ is closed in $B_\infty$.
\end{itemize}
Conversely, if  $\lim S_nf(x)$ exists  a.e. for all $f\in L^\infty(X, \mathcal{B}, \mu)$, then $S^*f(x) <\infty$ a.e. for every $f\in L^\infty(X, \mathcal{B}, \mu)$ and $S^*: B_\infty \to L^0(X, \mathcal{B},\mu)$ is continuous at $0$.
}
\medskip

Compared with the Banach principle, the Bellow-Jones principle requires the extra condition: the continuity at $0$ of $S^*$, meaning that $S^*f_n$ tends to $0$ in probability when $f_n \in B_\infty$ tends to $0$ in probability. If we apply the Banach principle to $\mathbb{B}=L^\infty(X, \mathcal{B}, \mu)$, the conclusion (B1) say that $S^*f_n$ tends to $0$ in probability when $\|f_n\|_\infty \to 0$.  The last assertion is weaker than the extra condition,  because $\|f_n\|_\infty \to 0$ is stronger than that $f_n$ tends to $0$ in probability. 

\begin{thm}[Characterization of $L^\infty$-Khintchin sequence]\label{thm:NSCinf}
    A sequence $(\lambda_n)$ of integers is a $L^\infty$-Khintchin sequence if and only if
    $T^*f(x)<\infty$ $\mathbf{m}$-a.e. for every $f\in L^\infty(\mathbb{T})$ and $T^*: B_\infty \to B_\infty$ is continuous at $0$, where $B_\infty$ is equipped with the topology of convergence in probability. 
\end{thm}
\begin{proof}
The necessity follows directly from the Bellow-Jones principle.

To prove the sufficiency, 
let us apply Bellow-Jones principle to our operators $(T_n)$ for the space $L^\infty(\mathbb{T}, \mathcal{B}, \mathbf{m})$.   Observe that continuous functions in $B_\infty$ are dense in $B_\infty$
and  the limit $\lim_{n\to \infty}T_nf(x)$ exists a.e. for every continuous function $f$.
So, for every $f\in L^\infty(\mathbb{T})$, $\lim_{n\to \infty}T_nf(x)$ exists a.e. To confirm that
the limit is equal to $\int_\mathbb{T}f(x)dx$ we can still apply Theorem \ref{thm:L^p-conv}, because 
$f\in L^p(\mathbb{T})$ for each $1\le p<\infty$. 
\end{proof}

\medskip

The above characterizations remain valid  for Khintchin sequences of group endomorphisms on an arbitrary compact abelian group that we discuss  in the next section.

    \subsection{Khintchin property depends on the ordering}
	
	First we observe that the ordering matters when the Khintchin class is in question.  Just look at $\mathbb{N}$ (naturally and increasingly ordered). Marstrand proved 
	that $\mathcal{K}_{\mathbb{N}}\not\supset L^\infty(\mathbb{T})$. However, reordering $\mathbb{N}$ can make $\mathbb{N}$ to be a Khintchin sequence.
	\begin{thm}
	     There is a reordering $\{r_n\}$ of $\mathbb{N}$ such that $\mathcal{K}_{\{r_n\}}=L^1(\mathbb{T})$.
	\end{thm}
	
	\begin{proof} Let $P$ be the increasing sequence of   the powers $2^n$ ($0\le n<\infty$). 
	Let $P'$ be the increasing sequence of   powers $2^{3^{m}}$ ($0\le m<\infty$), a very "small" subsequence of  $P$. 
	Enumerate $\mathbb{N}\setminus (P\setminus P')=(\mathbb{N}\setminus P)\sqcup P'$ as $\{b_0, b_1, \cdots, b_m, \cdots\} $ in an increasing order. Then  we replace $2^{3^{m}}$ in $P$ by $b_m$
	to get a sequence $P^*=\{r_n\}$, a reordering of $\mathbb{N}$.
	Assume that $f\in L^1(\mathbb{T})$. By Birkhoff's theorem, we have
	\begin{equation}\label{eq:Bk}
	    a.e. \quad \lim_{n\to \infty} \frac{1}{n}\sum_{k=0}^{n-1} f(2^k x) =\int f d\mathbf{m} .
	\end{equation}
    This is a nice equality. When $k=3^m$, the term $f(2^{3^m}x)$ will be kicked out. We will show that the equality remains true, see \eqref{eq:Bk2} below.  When the term $f(2^{3^m}x)$ is kicked out, we replace it by 
    $f(b_m x)$. We will also show that this modification doesn't affect the equality either, see \eqref{eq:Bk3}.

	As $ f(2^k x) = o(k)$ holds almost everywhere (a consequence of \eqref{eq:Bk}), we have $ f(2^{3^{m}} x) = o(3^{m})$ so that
	\begin{equation}\label{eq:Bk2}
	    a.e. \quad \lim_{n\to \infty} \frac{1}{n}\sum_{m: {3^m}<n} f(2^{3^{m}} x) =0.
	\end{equation}
	On the other hand, as $\sum n^{-2}\int |f(nx)| dx<\infty$, we have $f(nx) =O(n^2)$ a.e. so that
	\begin{equation}\label{eq:Bk3}
	    a.e. \quad \lim_{n\to \infty} \frac{1}{n}\sum_{m: 3^{m}<n} f(b_m x) =0,
	\end{equation}
	because $b_m=O(m^2)$  and 
	$
	 \sum_{k: 3^{m}<n} b_m^2  
	$
	is bounded up to a constant by   \[\sum_{k: 3^{m}<n} m^4=o(\log^5 n).\] 
	Combining \eqref{eq:Bk},  \eqref{eq:Bk2} and  \eqref{eq:Bk3} leads to 
	\begin{equation*}
	    a.e. \quad \lim_{n\to \infty} \frac{1}{n}\sum_{k=0}^{n-1} f(r_k x) =\int f d\mathbf{m} .
	\end{equation*}
	\end{proof}

	 In what follows, sequences of integers are always assumed increasingly ordered.
	

\subsection{Union and intersection of Khintchin sequences}

Consider two increasing sequences of integers $A=(a_n)$ and $B=(b_n)$. By $A \vee B$ we denote the increasing sequence obtained by re-ordering 
the integers in the union $A\cup B$.  Notice that duplicate integers are counted once and we always order the integers in $A \vee B$ in the increasing order. We are going to show that if $(a_n)$ and $(b_n)$ are Khintchin sequences, so is $(a_n)\vee (b_n)$. Subsequences of a Khintchin sequence are also considered.

Concerning the subsequences, we need the notion of relative density.  
  Let $\Lambda \subset \mathbb{N}$ and let $E\subset \Lambda$ be a subset of $\Lambda$. We define the  lower density of $E$ in $\Lambda$ by
$$
  \underline{d}_\Lambda(E) = \varliminf_{N\to \infty} \frac{|E\cap [1, N]|}{|\Lambda \cap [1, N]|}.
$$
The  upper density $\overline{d}_\Lambda(E) $ of $E$ in $\Lambda$  is similarly defined, by limsup.  If $\underline{d}_\Lambda(E) =\overline{d}_\Lambda(E)$,
we denote the common value by $d_\Lambda(E)$, called the density of $E$ in $\Lambda$.

When $\Lambda=\mathbb{N}$, we will talk about lower density $\underline{d}(E):=\underline{d}_{\mathbb{N}}(E)$ etc.

\begin{thm} \label{thm:union} Let $A=(a_n)$ and $B=(b_n)$ be two  $L^p$-Khintchin sequences for some $1\le p\le \infty$.  Let $A'\subset A$ be a subsequence of $A$.
Then \\
\indent {\rm (a)} $A\vee B$ is  a $L^p$-Khintchin sequence.\\
\indent {\rm (b)} $A'$ is a $L^p$-Khintchin sequence if $\underline{d}_A(A')>0$.\\
\indent {\rm (c)} $A\cap B$ is  a $L^p$-Khintchin sequence under the condition $\underline{d}_{A}(A\cap B)>0$.

\end{thm} 

\begin{proof} 
We first introduce the following notations. 
Let $f\in L^p(\mathbb{T})$. Denote 
 \begin{eqnarray*}
        A_nf(x) &=& \frac{1}{|A\cap [1,n]|}  \sum_{a \in A\cap [1, n]} f(ax),\\
        B_nf(x) &=&   \frac{1}{|B\cap [1,n]|}  \sum_{a \in B\cap [1, n]} f(ax), \\
        U_nf(x) & = & \frac{1}{|(A\cup B)\cap [1, n]|}\sum_{a \in (A\cup B)\cap [1, n]} f(ax).
    \end{eqnarray*}
    Let 
    $$
      A^*f(x):=\sup_{n\ge 1}|A_nf(x)|,  \quad B^*f(x):=\sup_{n\ge 1} |B_nf(x)|, 
      \quad U^*f(x):=\sup_{n\ge 1}|U_nf(x)|. 
    $$

    {\em Proof of (a) when $1\le p<\infty$}. 
     By Theorem~\ref{thm:NSCp}, it suffices to prove that
\[
U^{*}f(x) < \infty \quad \text{a.e. for every } f \in L^{p}(\mathbb{T}).
\]
Since the maximal operator is sublinear, we may assume without loss of generality that
\( f(x) \ge 0 \) almost everywhere. 
The key point is the observation that 
    \begin{equation*}\label{eq:Q}
        \sum_{a \in (A\cup B)\cap [1, n]} f(ax) =   \sum_{a \in A\cap [1, n]} f(ax)  +  \sum_{a \in B\cap [1, n]} f(ax)  -  \sum_{a \in (A\cap B) \cap [1, n]} f(ax).
    \end{equation*}

       As $f$ is nonnegative, we get   \begin{equation*}\label{eq:QABC}
     U_nf(x)\le  \frac{|A\cap [1, n]|}{|(A\cup B)\cap [1, n]|} \cdot A_nf(x) +  \frac{|B\cap [1, n]|}{|(A\cup B)\cap [1, n]|} \cdot B_nf(x).
       \end{equation*}
Then
    \begin{equation*}\label{eq:union}
      \forall\  0\le f\in L^p(\mathbb{T}), \quad  \ \ \  U^*f(x) \le A^*f(x) + B^*f(x) <\infty \ \ a.e. 
     \end{equation*}
     because  $A$ and $B$ are $L^p$-Khintchin sequences, which implies by Theorem \ref{thm:NSCp} 
 \begin{equation*}\label{eq:AB}
      \forall\  0\le f\in L^p(\mathbb{T}), \quad   A^*f(x) <\infty\  a.e \ \text{and} \ \ B^*f(x) <\infty \ a.e.
          \end{equation*}
          
{\em Proof of (b) when $1\le p<\infty$}. It suffices to notice that for nonnegative $f\in L^p(\mathbb{T})$
$$
\frac{1}{|A'\cap[1,n]|}\sum_{a\in A'\cap[1,n]} f(x) \le \frac{|A\cap [1,n]|}{|A'\cap [1,n]|}\cdot \frac{1}{|A\cap[1,n]|}\sum_{A\cap [1,n]} f(ax).
$$
By the hypothesis and Theorem \ref{thm:NSCp}, we have $A^*f(x)<\infty$ a.e., 
so that
$$
\sup_{n}\frac{1}{|A'\cap[1,n]|}\sum_{a\in A'\cap[1,n]} f(x) \le C A^*f(x)<\infty  \ \ \ a.e.,
$$
for some constant $C>0$ depending on the relative density $\underline{d}_A(A')>0$. 
The term on the left-hand side is nothing but the maximal function associated to the subsequence $A'$.
We conclude by using Theorem \ref{thm:NSCp}.

{\em Proof of (a)  when $p=\infty$}.
The above estimate $U^*f(x)\le A^*f(x) +B^*f(x)$ remains true for non-negative functions $f\in L^\infty$. Then for complex valued functions $f\in L^\infty$ we have
$$
U^*|f|(x) \le 4[A^*|f|(x) +B^*|f|(x)].
$$
By the hypothesis and Theorem \ref{thm:NSCinf},  we have $A^*|f|(x)<\infty$ a.e. and $B^*|f|(x)<\infty$ a.e. for $f\in B_\infty$, so that  $U^*f(x)<\infty$ a.e.  By the hypothesis and Theorem \ref{thm:NSCinf},  $A^*: B_\infty \to B_\infty$ and $B^*: B_\infty \to B_\infty$ are continuous at $0$. Then, so is $U^*: B_\infty \to B_\infty$.  We conclude by using  Theorem \ref{thm:NSCinf}. 

{\em Proof of (b) when $p=\infty$}. The argument is similar to the proof of (a) for $p=\infty$.

  {\em Proof of (c)}. It is a consequence of (b), because $A\cap B$ is a subsequence of $A$.
  \end{proof}

    The following is a direct corollary of  Theorem \ref{thm:union} and Marstrand's result that $\mathbb{N}$ is not a $L^\infty$-Khintchin sequence.
    
    \begin{cor} If $\mathbb{N}$ is decomposed into two disjoint sets $\Lambda_1$ and $\Lambda_2$, then one of  the sets $\Lambda_1$ and $\Lambda_2$ is
   not a $L^\infty$-Khintchin sequence (the sequence consisting of the integers in the set  arranged in increasing order).
    \end{cor}
    
    The following corollary follows from the fact that  $\{2^n\}$ and $\{3^m\}$ are $L^1$-Khintchin sequences.
    
     \begin{cor} $\{2^n\}\vee\{3^m\}= \{1,2,3,4,8,9,16,27, 32, 64, 81, \cdots\}$ is a $L^1$-Khintchin sequence.
      \end{cor}

      The following corollary follows from the fact that  $\{2^n 3^m: n\ge 0, m\ge 0\}$ and   $\{3^n 5^m: n\ge 0, m\ge 0\}$ are $L^1$-Khintchin sequences.

       \begin{cor} $\{2^n 3^m: n\ge 0, m\ge 0\}\vee \{3^n 5^m: n\ge 0, m\ge 0\}$ is a $L^1$-Khintchin sequence.
      \end{cor}

\subsection{Construction by multiplication and suggested questions}
 There are different ways to construct more or less regular sequences of integers and we investigate their Khintchin classes.  Let us introduce the following way. 
 First recall that every integer $a$ with $|a|\ge 2$ defines an expanding surjective endomorphism $\tau_a(x) = a x \mod 1$ on $\mathbb{T}$ which preserves the Lebesgue measure. 
 Consider the product space $\Omega:=(\Z\setminus \{-1, 0, 1\})^\N$.
 Any sequence $(\omega_n) \in \Omega$ produces a (lacunary) sequence of integers $(\omega_n\cdots \omega_2\omega_1)_{n\ge 1}$.  The following question is natural.

 \begin{Question}  \label{Q:1}
 {\em 
 What is the Khintchin class of this sequence $(\omega_n\cdots \omega_2\omega_1)_{n\ge 1}$? Under what condition is  $(\omega_n\cdots \omega_2\omega_1)_{n\ge 1}$ a Khintchin sequence?
 }
 \end{Question}
 
  The Khintchin class of $(\omega_n\cdots \omega_2\omega_1)_{n\ge 1}$
 depends on $\omega \in \Omega$.  
 The  examples  mentioned in Introduction are all special cases of the above construction.  Raikov example $(2^n)$ corresponds to the choice $\omega_n=2$ for all $n\ge 1$ and 
 the Rosenblatt's example $(2^{2^n})$ corresponds to the choice $\omega_n=2^{2^{n-1}}$ for all $n\ge 1$.  If we choose $\omega_n =2^{2n-1}$ for all $n\ge 1$,
 we get $\omega_n\cdots \omega_2\omega_1=2^{n^2}$ and the corresponding sequence $(2^{n^2})$ was investigated by Bourgain,  Buczolich and Mauldin.
 \medskip
 
 The space $\Omega$ is huge (not compact under the product topology). It seems unrealistic to answer Question \ref{Q:1} for every $\omega$ in it. 
   We can consider a space $\Omega_0$ smaller than $\Omega$.  Then Question \ref{Q:1} with $\omega\in\Omega_0$ would be simpler.  For example, let us fix two integers, say $2$ and $3$, and consider $\Omega_0:=\{2,3\}^\mathbb{N}$.

One of sub-questions of Question \ref{Q:1} that we ask is stated as follows. 
Let us take $2$ or $3$ independently to get a random sequence $(\omega_n)$ of $2$ or $3$. More precisely,
we assume that $\omega_n$'s are independent and $P(\omega_n =2)=p=1-P(\omega_n =3)$ ($0<p<1$ being fixed). Then we say  that   $(\omega_n)$ is
a {\em $p$-Bernoulli sequence of $2$ or $3$}. 

\begin{Question}\label{Q:2}
Let 
$(\omega_n)$ be
a \em $p$-Bernoulli sequence of $2$ or $3$.
Is  $(\omega_n\cdots \omega_2\omega_1)_{n\ge 1}$ almost surely a Khintchin sequence?
\end{Question}

We could conjecture that $(\omega_n\cdots \omega_2\omega_1)_{n\ge 1}$ is almost surely a $L^p$-Khintchin sequence for some $p\ge 1$.
We are not able to prove the conjecture, but we can prove the following result, which supports the conjecture but is still far from the conjectured result.

\begin{thm} \label{thm:quasiKP} If $(\omega_n)$ is
a $p$-Bernoulli sequence of $2$ or $3$ for some $0<p<1$, then for any $f\in L^1(\mathbb{T})$, we have $f\in \mathcal{K}_{\{\omega_n\cdots \omega_1\}}$ a.s. 

\end{thm} 

Theorem \ref{thm:quasiKP} is a very special case of a result that we shall prove (see Theorem \ref{thm:WKS}).
Actually a very large class of candidates of random  Khintchin sequences can be constructed by using 	Anzai's skew product construction, even on any compact abelian group. 
\medskip

The space $\Omega_0 = \{2,3\}^{\mathbb{N}}$ is a compact metric space.
The shift map $\sigma : \Omega_0 \to \Omega_0$ is defined by
$(\sigma \omega)_n = \omega_{n+1}$, and  is continuous.
 Thus we get a topological dynamical system 
$(\Omega_0, \sigma)$. 
The space $\Omega_0$ is still big. 
We can consider subspaces of $\Omega_0$, like subshifts, i.e. closed $\sigma$-invariant subspace $\Omega_{00}$ (i.e. $\sigma(\Omega_{00})\subset \Omega_{00}$).  There are many such subsystems
$(\Omega_{00}, \sigma)$. They can be minimal (i.e., every orbit $(\sigma^n \omega)_{n\ge 0}$ is dense in $\Omega_{00}$) or uniquely ergodic (i.e., there is a unique $\sigma$-invariant probability measure on $\Omega_{00}$). 

\begin{Question}\label{qest:MUE}
{\em 
Let $\Omega_{00}$ be a subshift of $\Omega_0$ which is minimal and unique ergodic.  Is\\  $(\omega_n\cdots \omega_2\omega_1)_{n\ge 1}$ a $L^p$-Khintchin sequence
for every $\omega \in \Omega_{00}$?
}
\end{Question}

Yes, we will prove that it is the case for the Thue-Morse subshift (see Theorem \ref{thm:TM}). This is the first such example that we know. It has since been proved in \cite{FFL} that it is true for 
 the class of subshifts associated to balanced primitive substitutive sequences. For example, it is the case for the Fibonacci  sequence $(\omega_n)\in \{2,3\}^\infty$   of $2$ and $3$.   

Recall that $\Omega_0$ is a complete metric (compact) space. Let us consider the set of sequences in $\Omega_0$, denoted by $K_{2, 3}$, which produce Khintchin sequences. 

\begin{Question}\label{quest:K23}
{\em 
Do we have $K_{2,3} \not= \Omega_0$ (i.e. the existence of sequence  of $2$ and $3$ such that their first products don't form a Khintchin sequence)?
Is $K_{2,3}$ a set   of second Baire category in $\Omega_0$?
}

\end{Question}

The first question in Question~\ref{quest:K23} has an affirmative answer in \cite{FFL}, where the authors have constructed non-Khintchin sequences in $\Omega_0=\{2,3\}^\infty$ by using a Rokhlin lemma due to Avila and Candela \cite{AC2016}.  But we have no answer to the second question.

Similarly we denote by  $K_{2,3,4, \cdots}$ the set of sequences in $\{2,3, 4, \cdots\}^\infty$ (any integer larger than $2$ can be chosen) which produce by multiplication Khintchin sequences. 
In this case, we know that $K_{2,3,4, \cdots}$ is a proper subset  of $\{2,3, 4,\cdots\}^\infty$.   Rosenblatt's counter-example is not in $K_{2,3,4, \cdots}$. 

\begin{Question}
Is $K_{2,3, 4,  \cdots}$ a second Baire category set in $\{2,3, 4, \cdots\}^\infty$ ?
\end{Question}

Let us add the following question to our list of questions.  

\begin{Question}
	Assume that $(a_n)$ is a Khintchin sequence. Is the translated sequence $(a_n+a)$ (with $a\in \mathbb{Z}\setminus\{0\}$) still 
	 a Khintchin sequence? In particular, is $(2^n +1)$
     a Khintchin sequence?
	\end{Question}

\medskip

	\subsection{Thue-Morse sequence is a Khintchin sequence} 
	
	Treat $2$ and $3$ as letters in the alphabet $\{2,3\}$ and then consider the substitution $$2 \mapsto 23,  \quad 3\mapsto 32$$
	 which produce the Thue-Morse sequence $(t_n)$:
$$
    2, 3, 3, 2, 3,2,2, 3, 3,2,2,3,2,3,3,2, \cdots
$$   
(cf. \cite{Queffelec} for substitutions). Then consider the lacunary sequence $(t_n^*)$ defined by 
the products $$t_n^* = t_n\cdots t_2 t_1.$$
We call $(t_n^*)$ the {\em Thue-Morse product sequence} of $2$ and $3$. 
The first terms of $(t_n^*)_{n\ge 1}$ are:
$$
     2^1\cdot 3^0, 
     {2\cdot3}, 2\cdot3^2,  
     {2^2\cdot3^2}, 2^2\cdot3^3, 
     {2^3\cdot3^3}, 2^4\cdot 3^3,  
     {2^4\cdot3^4}, 2^4\cdot3^5,  
     {2^5\cdot3^5},  2^6\cdot3^5,  
     {2^6\cdot3^6},  2^7\cdot3^6,
      \cdots
$$ 
or 
\begin{equation}\label{eq:TM}
     \textcolor{blue}{2}\cdot 6^0, \textcolor{red}{6},  \textcolor{green}{3}\cdot 6; \ \  \textcolor{red}{6^2},  \textcolor{green}{3} \cdot6^2; \ \ \textcolor{red}{6^3},  \textcolor{blue}{2} \cdot 6^3; \ \ \textcolor{red}{6^4},  \textcolor{green}{3} \cdot6^4; \ \ \textcolor{red}{6^5},   \textcolor{blue}{2} \cdot6^5; \ \ \textcolor{red}{6^6},    \textcolor{blue}{2} \cdot6^6; \ \
     \cdots
\end{equation}
We see that there are three types of terms in $(t_n^*)_{n\ge 1}$:   $a\cdot 6^n$ with $a=1,2,3$. Define 
$$
     \mathcal{U}_a= \{ n\ge 0: a\cdot 6^n =t_m^* \ {\rm for \ some} \ m\}, \qquad (a=1,2,3).
$$
Clearly the union of $6^{\mathcal{U}_1}$, $2\cdot 6^{\mathcal{U}_2}$ and $3 \cdot 6^{\mathcal{U}_3}$  is the set of all $\lambda_n$ with $n\ge 1$.

\begin{lem} 
The densities of $ \mathcal{U}_1, \mathcal{U}_2$ and $ \mathcal{U}_3$
are respectively equal to $d_1=\frac{1}{2}, d_2=\frac{1}{4}$ and $d_3=\frac{1}{4}$.
\end{lem}
\begin{proof} 
The result is well known for experts on the subject. We give a proof for the reader's convenience. As the sequence $(t_n)$ is a sequence of $23$ and $32$, clearly $\mathcal{U}_1$ has the density
$\frac{1}{2}$. Notice that  $n \in \mathcal{U}_2$ iff $n=2m+1$ and $t_{2m+1}=2$, and  $n \in \mathcal{U}_3$ iff $n=2m+1$ and $t_{2m+1}=3$. 
However it is known that $(t_{2m+1})_{m\ge 0}$ is a Thue-Morse sequence too and it is also known that both frequencies of $2$ and $3$ in $(t_n)$ are equal to $\frac{1}{2}$. 
So, both $\mathcal{U}_2$ and $\mathcal{U}_3$ has $\frac{1}{4}$ as density.
\end{proof}

We are now ready to state and prove the following result.
\begin{thm}\label{thm:TM}  The Thue-Morse product sequence $(t_n^*)_{n\ge 1}$ is a $L^1$-Khintchin sequence.
\end{thm}

\begin{proof} For $f \in L^1(\mathbb{T})$ and $N\ge 1$, we have
$$
     \sum_{n=1}^N f(t_n^* x)
     =   \sum_{k \in \mathcal{U}_1, 6^k\le N}^N f(6^k x) + \sum_{k \in \mathcal{U}_2, 2\cdot 6^k\le N}^N f(2\cdot 6^k x) + \sum_{k \in \mathcal{U}_3, 3\cdot 6^k\le N}^N f(3\cdot 6^k x).
$$
By Theorem \ref{thm:union} (b) and the fact that the Lebesgue measure is ergodic with respect to $x \mapsto 6x \mod 1$, we have
$$
    a.e. \ \ \ \lim_{N\to \infty} \frac{1}{N} \sum_{k \in \mathcal{U}_a, a\cdot 6^k\le N}^N f(a\cdot 6^k x) = d_a \int f(ay) dy = d_a\int f(z) dz
$$
for $a=1,2$ and $3$. It follows that
$$
a.e. \ \ \ \lim_{N\to \infty} \frac{1}{N}  \sum_{n=1}^N f(t_n^* x) = (d_1+d_2+d_3) \int f(z) dz = \int f d\mathbf{m} .
$$
\end{proof}

Here the proof is based on the specific structure of the Thue-Morse sequence. The key point is that we have only powers of $6$ in \eqref{eq:TM}.  This argument does not work for the Fibonacci sequence for which there is no basis like $6$ for the Thue-Morse sequence. However, 
the following result is now proved in \cite{FFL}.
\begin{thm}[\cite{FFL}]
    Any $C$-balanced primitive substitutive sequence produces a Khintchin sequence in the way we have introduced above. 
\end{thm}

Thus Question \ref{qest:MUE} is answered affirmatively for the subshift defined by a $C$-balanced sequence.

Recall that a {\em substitution} $\sigma$ is defined as a non-erasing morphism of the free monoid over a given finite alphabet $A$ which admits an infinite word fixed by 
$\sigma$; it is {\em primitive} if there exists an integer
$n$ such that the image by the $n$-th iterate $\sigma^n$ of $\sigma$ of any letter of the alphabet contains all the letters of the alphabet. 

An infinite word $u$ over the alphabet $A$ is said to be {\em $C$-balanced} if for any
pair of factors $(w,w')$ of $u$ of the same length and for any letter $a$ of the alphabet, the difference 
$||w|_a-|w'|_a|$
 between the numbers of occurrences of the letter $a$ in those two factors is bounded by $C$ in 
absolute value for some constant $C\ge 1$. When $C=1$, the word $u$ is usually said to be {\em balanced}. The balance function $B(n)$ of an infinite word $u$ is then defined as 
$$
B(n) =\max_{a\in A} \max_{|w|=|w'|=n} ||w|_a-|w'|_a|
$$
(the supremum being taken on the factors $w, w'$ of $u$), where $|w|_a$ denotes the number of occurrences of the letter $a$ in the word $w$.
Adamczewski \cite{Adamczewski2003} gave a precise description in terms of the second largest eigenvalue of the incidence matrix of the substitution $\sigma$ of the behaviour of the balance function. In particular, an infinite word fixed by a primitive substitution which is  C-balanced for some constant C if and only if the second largest eigenvalue of the incidence matrix is  of modulus smaller than 1 (this is the Pisot case), or it is a simple eigenvalue of modulus one with some extra combinatorial condition which does not only depend on the incidence matrix, but also on the order of occurrence of the letters (in the images of the letters by the substitution). 

Our Fibonacci sequence is defined by 
\[2\mapsto 23, \quad 3\mapsto 2.\]
The first terms of the Fibonacci sequence are 
\[2,3,2,2,3,2,3,2,2,3,2,2,3,2,3,2,2,3,\cdots.\]

    \subsection{Regular functions and Khintchin classes}
    A sequence of integers $(\lambda_n)$ is Hadamard lacunary if $\inf_k \frac{\lambda_{k+1}}{\lambda_k}\ge q>1$ for some $q$.  Let us state the following result due to Erd\"{o}s which states that a function of some regularity belongs to the Khintchin class $\mathcal{K}_{\{\lambda_n\}}$.

    \begin{thm}[Erd\"{o}s \cite{Erdos1949}]
        Let $(\lambda_n)$ be a Hadamard lacunary sequence of integers. Let $f\in L^2(\mathbb{T})$. Then $f \in 
        \mathcal{K}_{\{\lambda_n\}}$ if there exist constants $A>0$ and $\alpha>0$ such that for large $N$ we have
        \begin{equation}\label{eq:Erdos}
        \sum_{|n|\ge N} |\widehat{f}(n)|^2 \le \frac{A}{(\log \log N)^\alpha}.
        \end{equation}
    \end{thm}
    
As Rosenblatt's example shows, in general case, for $f$ to belong to the Khintchin class, some regularity of $f$ is necessary. Erd\"{o}s' condition \eqref{eq:Erdos} is such a sufficient condition. 

    It is interesting to compare the above Erd\"{o}s result to a corollary of a result due to Cuny and Fan \cite{CunyFan}.

\begin{thm}[\cite{CunyFan}]\label{thm:CF}
        Let $(\lambda_n)$ be a Hadamard lacunary sequence of integers. Let $f\in L^2(\mathbb{T})$ with $\int_\mathbb{T}f(x) dx=0$. Suppose that there exists  $\epsilon>0$ such that
\begin{equation}\label{eq:CunyFan}
\sum_{|n|>N} |\widehat{f}(n)|^2 = O\left(\frac{1}{\log^{1+\epsilon} N}\right).
\end{equation}
Then for any square-summable sequence $(a_n)$, the series
$\sum a_n f(\lambda_n x)$ converges almost everywhere. 
    \end{thm}
    \begin{proof}
        A result in \cite{CunyFan} asserts that the conclusion of Theorem \ref{thm:CF} holds true if $$
	\omega_{f,2}(\delta)  = O\left(\frac{1}{\log^{\alpha/2} \frac{1}{\delta}}\right)
	$$
	for some $\alpha >1$, where $\omega_{f,2}(\cdot)$
    is the $L^2$-modulus of continuity of $f$ defined by
    $$
    \omega_{f, 2}(\delta)= \sup_{|h|\le \delta}\|f(\cdot +h)-f(\cdot)\|_2.
    $$
    
    We are going to check this condition by using the hypothesis made on the Fourier coefficients of $f$. Indeed, 
	$$
	   \|f(\cdot + h)-f(\cdot)\|_2^2 = \sum_{n} |\widehat{f}(n)|^2 |e^{2\pi i nh} -1|^2= 4  \sum_{n} |\widehat{f}(n)|^2 |\sin \pi n h |^2.
	$$
	Let  
	$
	      R(N) = \sum_{|n|>N} |\widehat{f}(n)|^2.
	$
	Take the integer $N =[h^{-1/2}]$ so that $Nh\le h^{1/2}$ and
	$$
	    \|f(\cdot + h)-f(\cdot)\|_2^2  \le 4  \sum_{|n|\le N}  |\widehat{f}(n)|^2 (nh)^2 + 4 R(N) \le 4 \|f\|_2^2 h + 4 R(h^{-1/2}).
	$$
	As $R(N) = O(1/\log^\alpha N)$, the main term in the above estimate is 
$R(h^{-1/2})$. Then 
	 we obtain the desired estimate.
      \end{proof}

	Gaposhkin \cite{Gaposhkin1967} produced an example showing that the condition \eqref{eq:CunyFan} is optimal, namely $\epsilon$ cannot be reduced to $0$.

\section{Khintchin sequences of group endomorphisms}\label{sect:groups}

In this section, we generalize the notion of  Khintchin sequence to arbitrary compact abelian groups. Endomorphisms on such groups need not commute, which makes the problem more challenging.

\subsection{Epimorphisms on compact abelian groups}
	
	Consider  a compact abelian  group $G$ with dual group $\widehat{G}$ and Haar measure $dx$, sometimes denoted $\mathbf{m}$. 
	Let $C(G)$ be the Banach space of continuous function on $G$, equipped with the supremum norm $\|\cdot\|_\infty$. The dual space $C(C)^*$
	of $C(G)$
	represents the space  $M(G)$ of finite Borel measures on $G$. 
	A sequence of points $(x_n)_{n\ge 0}$ on $G$ is said to be {\em equidistributed} or {\em uniformly distributed} if 
	the probability measures
	$n^{-1}\sum_{k=0}^{n-1}\delta_{x_k} \in M(G)$ converges in weak$^*$-topology to the Haar measure, namely
	for every continuous function $f\in C(G)$ we have
	\begin{equation}\label{eq:ud-def}
	    \lim_{n\to \infty} \frac{1}{n}\sum_{k=0}^{n-1} f(x_n) =\int_G f  d\mathbf{m} .
	\end{equation}
	An equivalent statement (Weyl's criterion) is that for any $\gamma \in \widehat{G}\setminus\{1\}$ we have  
	$$
	    \lim_{n\to \infty} \frac{1}{n}\sum_{k=0}^{n-1} \gamma(x_n) =0.
	$$ 
    We refer to \cite{Kuipers1974} for more information about uniform distributions on compact groups, even non abelian.  
	
	We denote by ${\rm End}(G)$ the semigroup of all continuous endomorphisms on $G$, by  ${\rm Epi}(G)$
		 the semigroup of all continuous surjective endomorphisms on $G$ (also called epimorphisms), and by ${\rm Aut}(G)$ the group of all continuous automorphisms of $G$.  
		 Any epimorphism preserves the Haar measure. Let us recall the simple proof of this fact. Let $\tau$ be an epimorphism on $G$
	and $B$ a Borel set of $G$. For any $y\in G$ we have some $y'\in G$ such that $\tau (y')=y$ so that 
	\begin{eqnarray*}
	\mathbf{m}\circ \tau^{-1}(B - y ) & =&  \int 1_{B}(\tau (x) +y) d\mathbf{m} (x)= \int 1_{B}(\tau (x) +\tau(y') d\mathbf{m} (x) \\
	&=& \int 1_{B}(\tau (x+y') d\mathbf{m} (x)=\int 1_{B}(\tau (x) d\mathbf{m} (x)= \mathbf{m}\circ \tau^{-1}(B).
	\end{eqnarray*}
	
	 Let us look at two typical examples of compact abelian groups. 
	 \medskip
		 
		{\bf Example 1.}  It is known that for $G=\mathbb{T}^d$, the d-dimensional torus, we have
		 $$
		     {\rm End}(\mathbb{T}^d) \simeq M_d(\mathbb{Z}), \quad  {\rm Epi}(\mathbb{T}^d) \simeq GL_d(\mathbb{Z}), 
		     \quad  {\rm Aut}(\mathbb{T}^d) \simeq SL_d(\mathbb{Z}),
		 $$
		 where $M_d(\mathbb{Z})$ is the set of integral matrices (i.e. matrices with entries in $\mathbb{Z}$), $GL_d(\mathbb{Z})$ is the set of integral matrices $M$ such that  $\det M\not=0$
		 and $SL_d(\mathbb{Z})$ is the set of integral matrices $M$ such that  $\det M=\pm 1$ (cf. \cite{Walters1982},p.15).  So, an endomorphism on $\mathbb{T}^d$ is of the form
		 $x \mapsto Ax \mod \mathbb{Z}^d$ for some matrix $A\in M_d(\mathbb{Z})$. Those endomorphisms preserving the Haar measure are the surjective endomorphisms.
		 Some of surjective endomorphisms are ergodic with respect to the Haar measure, which are characterized by the fact that the corresponding matrix has no root of unit as eigenvalue
		 (cf. \cite{Walters1982}). 
		 \medskip
		 
		 {\bf Example 2.}   For the additive group $\mathbb{Z}_p$ of $p$-adic numbers, we have 
		  $$
		     {\rm End}(\mathbb{Z}_p) \simeq  \mathbb{Z}_p, \quad  {\rm Epi}(\mathbb{Z}_p) =  {\rm Aut}(\mathbb{Z}_p) \simeq \mathbb{U},
		 $$
		 where $\mathbb{U}=\{x\in \mathbb{Z}_p: |x|_p=1\}$ is the group of unit of $\mathbb{Z}_p$, where $\mathbb{Z}_p$ is  considered as the ring of $p$-adic integers.
		 So, an epimorphism (or equivalently automorphism) on $\mathbb{Z}_p$ is of the form $x\to a x$ with $a\in \mathbb{U}$. But, notice that no epimorphism on $\mathbb{Z}_p$ is ergodic with respect to the Haar measure because  $p^m \mathbb{Z}_p$
		 are all invariant sets of positive measure for $m\ge 1$. We refer to \cite{Koblitz1984} for the information about the group $\mathbb{Z}_p$. 
		 \medskip
		 
		 Return back to a general compact abelian group $G$. 
		 Any endomorphism $\tau\in {\rm End}(G)$ induces 
		 its adjoint, a continuous endomorphism $\tau^*: \widehat{G}\to \widehat{G}$ on $\widehat{G}$, namely $$
		 \tau^*(\gamma)= \gamma \circ \tau \ \ \ {\rm for} \ \ \gamma \in \widehat{G}.
		 $$ 
		 It is easy to see that $\tau^*$ is monic (resp. epic) iff  $\tau$ is epic (resp. monic). Here we use ``monic" as synonym of ``injective" and ``epic" as synonym of ``surjective". Our endomorphisms will all be assumed epic, so their adjoints will  all be monic.    
		 Since $\widehat{G}$ is discrete, the dynamical system $(\widehat{G}, \tau^*)$ is relatively simple and it serves as a tool to study the actions of epimorphisms on $G$. For example, it is well known that  an epimorphism  is ergodic (with respect to the Haar measure)	
		 iff the $\tau^*$-orbit of every character $\gamma\not=1$ is infinite (here $1$ denotes the identity element of the dual group for which the operation is assumed multiplicative).  
\medskip

\subsection{Concepts of Khintchin class and Khintchin  sequence}	
	Suppose we are given a sequence of epimorphisms  $\{\tau_n\}\subset {\rm Epi}(G)$ of a compact abelian group $G$. 
	 We will study the equidistribution (i.e. uniform distribution) on the group $G$ of the sequence $\{\tau_n x\}_{n\ge 1}$ for $x\in G$. 
	The first question is to find condition on $\{\tau_n\}$ such that  $\{\tau_n x\}_{n\ge 1}$  is uniformly distributed for almost all $x\in G$. 
	If  $\{\tau_n x\}_{n\ge 1}$  is uniformly distributed for almost all $x\in G$, we say that  $\{\tau_n\}$ is a {\em ud-sequence} of epimorphisms. 
	Following \cite{AHM} which studied {Khintchin classes} in $\mathbb{Z}$, for a ud-sequence $\{\tau_n\} \subset {\rm Epi}(G)$, we define its {\em Khintchin class} by
	$$
	     \mathcal{K}_{\{\tau_n\}}=\left\{f\in L^1(G): \lim_{n\to \infty} \frac{1}{n}\sum_{k=0}^{n-1} f(\tau_n x) =\int f d\mathbf{m} \ \ a.e. \right\}.
	$$
	As $\{\tau_n\}$ is a ud-sequence, we have $ \mathcal{K}_{\{\tau_n\}} \supset C(G)$. Also, if $B$ is a Borel set in $G$ such that $\mathbf{m}(\partial B)=0$,
	we have $1_B \in  \mathcal{K}_{\{\tau_n\}}$, where $\partial B$ denotes the boundary of $B$.  For general Borel set $B$, it may be not true that $1_B \in  \mathcal{K}_{\{\tau_n\}}$. 
	That is the case even when $G=\mathbb{T}$ (Marstrand \cite{Marstrand1970}), as we have already mentioned.
	
	When $	     \mathcal{K}_{\{\tau_n\}}=L^1(G)$, we say that $\{\tau_n\}$ is a {\em Khintchin sequence} {\color{red}  for $G$.}  
	When $	     \mathcal{K}_{\{\tau_n\}}$ contains a subspace $\mathbb{B}$ of $L^1(G)$, we say that $\{\tau_n\}$ is a $\mathbb{B}$-{\em Khintchin sequence}  {\color{red} for $G$. } 
	For example, we can talk about $L^p$-{\em Khintchin sequence}  for some $1\le p\le \infty$.
	
	\medskip

        It is well known that every sequence of distinct integers is a ud-sequence in $\mathbb{Z} =\widehat{T}$.  But this statement is not true for general compact abelian groups. For example,  there is no ud-sequence in $\widehat{\mathbb{Z}_p}$, 
	because the epimorphisms on $\mathbb{Z}_p$ are automorphisms, which preserves every open set
	$p^m \mathbb{Z}_p$ ($m\ge 1$).  As was mentioned before,   the set of natural integers $\mathbb{N}$, increasingly ordered,  is not a Khintchin sequence (Marstrand's theorem) and 
	the sequence$\{2^n\}$ is a Khintchin sequence (Birkhoff's theorem), as well as the Furstenberg sequence, the increasingly ordered sequence of the Furstenberg set
	$\{2^n3^m: n\ge 0, m\ge 0\}$.
	The situation in general groups, even on $\mathbb{T}^d$ with $d\ge 2$, is more complicated.

     \subsection{$L^p$-convergence}
	However,
	the $L^p$-convergence always holds for ud-sequences and for $L^p$-integrable functions. 
	
	\begin{thm} \label{thm:L^p-conv}
    Let $1\le p<\infty$. Let $(\tau_n)$ be a ud-sequence of epimorphisms on a compact abelian group $G$. Then for every $f\in L^p(G)$, we have
	$$
	     L^p\!-\!\lim_{N\to \infty} \frac{1}{N} \sum_{n=1}^N f(\tau_n x) =\int_G f d\mathbf{m}.
	$$
	\end{thm}
	\begin{proof} We assume that $\int f d\mathbf{m}=0$, without loss of generality. Let 
	$$
	A_Nf(x) :=\frac{1}{N} \sum_{n=1}^N f(\tau_n x). 
	$$
	It suffices to show that 
	\begin{equation}\label{eq:Lp-conv1}
	     \lim_{N\to \infty}\|A_Nf\|_p=0.
	\end{equation}
	Since $(\tau_n)$ is a ud-sequence, for any continuous function $f$, we have $A_Nf(x) \to 0$ almost everywhere, so that \eqref{eq:Lp-conv1} holds for continuous functions, by the Lebesgue dominated convergence theorem. 
	
	For a general function $f\in L^p(G)$ with $\int fd\mathbf{m}=0$, we  approximate $f$ by continuous functions as follows (there are other approximations). Take an open neighborhood $V_\epsilon$ of the identity element  $1_G$ with parameter
	$\epsilon>0$ such that $V_\epsilon$ decreases to $\{1_G\}$ when $\epsilon \downarrow 0$. Let 
	$$
	     f^{(\epsilon)}(x) := f * \eta_\epsilon (x)=\int f(y)\eta_\epsilon(x-y) dy, \quad {\rm where}\ \ \ 
	\eta_\epsilon(x) = \frac{1}{\mathbf{m}(V_\epsilon)}1_{V_\epsilon}(x).
	$$ 
	Observe that  for fixed $\epsilon>0$, $f^{(\epsilon)}(\cdot)$ is  a continuous function and that $\|f-f^{(\epsilon)}\|_p\to 0$ when  $\epsilon \to 0$.
	As 
	$$
	     \|A_Nf\|_p \le  \|A_Nf -  A_Nf^{(\epsilon)}\|_p  +\| A_N f^{(\epsilon)}\|_p,
	$$
	and as $\| A_Nf^{(\epsilon)}\|_p \to 0$ (\eqref{eq:Lp-conv1} being proved for continuous functions), 
	we have
	$$
	     \varlimsup_{N\to \infty}\|A_Nf\|_p \le   \varlimsup_{N\to \infty} \|A_Nf -  A_Nf^{(\epsilon)}\|_p.
	$$
	However, since $\tau_n$ preserves the Haar measure, we have
	\begin{eqnarray*}
	\|A_Nf -  A_N^{(\epsilon)}f\|_p
	&\le& \frac{1}{N} \sum_{n=1}^N \| f \circ \tau_n - f^{(\epsilon)}\circ \tau_n\|_p=\|f-f^{(\epsilon)}\|_p.
	\end{eqnarray*} 
	Therefore
	$$
	   	     \varlimsup_{N\to \infty}\|A_Nf\|_p \le \|f-f^{(\epsilon)}\|_p.
	$$
	Letting $\epsilon \to 0$ to finish the proof.
	\end{proof}
	
	\begin{cor} Let $(A_n)$ be a sequence of expanding integral $d\times d$ matrices. For any $f\in L^p(\mathbb{T}^d)$ ($1\le p<\infty$), the average
	    $$
	         \frac{1}{N}\sum_{n=1}^N f(A_n\cdots A_2A_1 x)
	    $$ 
	    converges in $L^p$-norm to $\int f d\mathbf{m}$.
	\end{cor}

    \begin{proof} It follows immediately from Theorem \ref{thm:L^p-conv}
    and Theorem \ref{prop:Exp-ud} proved below.
    The expanding property will be discussed later:
    a matrix is expanding iff its singular values are all larger than $1$ (cf. Lemma \ref{lem:exp}).     
    \end{proof}
	
	In particular, if $(\lambda_n)$ be an increasing sequence of natural numbers such that $\lambda_n |\lambda_{n+1}$ for all $n\ge 1$ (this divisibility is not really needed), then
	for any $f\in L^p(\mathbb{T})$, $N^{-1}\sum_{n=1}^N f(\lambda_n x)$ tends in $L^p$-norm to $\int f$.
	Rosenblatt's example shows that for general such $\{\lambda_n\}$,  the almost everywhere convergence of $N^{-1}\sum_{n=1}^N f(\lambda_n x)$ may be violated.
	
	Let us point out two works related to  Theorem  \ref{thm:L^p-conv}.  The first  one is  the   following  particular case of a result due to  Blum and Hanson \cite{BH1960}: if $\tau_n =\tau^{k_n}$
	for some strongly mixing automorphism $\tau$ of $G$ and for some strictly increasing sequence of positive integers $(k_n)$, then the conclusion of Theorem  \ref{thm:L^p-conv}  holds  for the sequence $(\tau_n)$.  A non-zero Gaussian integer $\mathbf{n}=a+i b$ with $a, b \in \mathbb{Z}$ defines an epimorphism of $\mathbb{T}^2$:
	$$
	 T_{\mathbf{n}} = \begin{pmatrix}
	           a & -b\\
	           b & a
	       \end{pmatrix} 
	$$
	In \cite{Easwaran2012}, the conclusion of Theorem  \ref{thm:L^p-conv} is proved for some sequence of such epimorphisms. 
	
	\subsection{The sequence of powers of an epimorphism}
	If $\tau\in {\rm Epi}(G)$ is nilpotent (i.e. $\tau^m =I$ for some $m\ge 1$), then $\{\tau^n\}$ is periodic and cannot be a ud-sequence.
	However, if $\tau$ is not nilpotent, we have the following classification.
	
	\begin{prop} Assume that $ \tau \in {\rm Epi}(G)$ is not nilpotent. Then $\{\tau^n\}$ is a ud-sequence of epimorphisms on $G$ iff $\tau$
	 is ergodic. In this case, we have $\mathcal{K}_{\{\tau^n\}}=L^1(G)$.
	 \end{prop}
	 \begin{proof}
	 The ergodicity of $\tau$ implies that  $\{\tau^n\}$ is a ud-sequence and $\mathcal{K}_{\{\tau^n\}}=L^1(G)$, by Birkhoff's ergodic theorem. On the other hand, the non-ergodicity of $\tau$ means that 
	for some  $\gamma \in \widehat{G}\setminus\{1\}$,  
    then $\tau^{*m}\gamma =\gamma$ for some $m\ge 1$ (we assume that $m$ is the least one). Then 
	the non -trivial trigonometric polynomial
	$$f=\sum_{k=0}^{m-1} \tau^{*k}\gamma = \sum_{k=0}^{m-1} \gamma\circ \tau^k$$
	is $\tau$-invariant. 
	Therefore, 
	$$
	 \forall x, \quad  \lim_{n\to \infty} \frac{1}{n}\sum_{k=0}^{n-1} f(\tau^n x) =f(x).
	$$ 
	 So,  $(\tau^n x)$ is not uniformly distributed for almost all $x$, because it is not true that $f(x)\not= \int fd\mathbf{m}$ a.e.
	\end{proof}

    Let us look at the special group $\mathbb{T}^d$.
	
	\begin{cor} Assume that $ \tau \in {\rm GL}_d(\mathbb{Z})$. Then $\{A^n\}$ is a ud-sequence of epimorphisms on $\mathbb{T}^d$ iff $A$
	has no root of unity as eigenvalue.
	In this case, we have $\mathcal{K}_{\{A^n\}}=L^1(\mathbb{T}^d)$.
	 \end{cor}
	
	The exceptional sets regarding the equidistributions of $\{A^n x\}$ were studied in \cite{Fan1991,Fan1993,Fan1995,Lesigne1998}.
	\medskip

	Let $\tau$ be an ergodic epimorphism on $G$ and let $\mathcal{U}=\{u_k\}$ be an increasing sequence of positive integers. Then we consider
	the sequence of epimorphisms $\tau^{\mathcal{U}}:=\{\tau_n\}$ defined by $\tau_n =\tau^{u_n}$.
	\medskip
	
	\begin{thm} \label{thm:positivedensity} Let $\tau$ be an ergodic epimorphism on $G$.
	The sequence $\tau^\mathcal{U}$ is a Khintchin sequence if $\mathcal{U}$ has a positive lower density, i.e.
	$$
	      \varliminf_{n\to \infty} \frac{|\mathcal{U}\cap [1, n]|}{n} >0.
	$$
	\end{thm}
	\begin{proof} 
    We cannot apply Theorem \ref{thm:union} (b) directly in this general compact-group setting,  because its proof relies on applying Theorem 3.1 to a subsequence of integers. Therefore, we give a self-contained argument.

    For trigonometric polynomials $f$, which are dense in $L^1(G)$, we have
	\begin{equation}\label{eq:tau-powers}
	     a.e. \quad \lim_{n\to \infty} \frac{1}{|\mathcal{U}\cap [1, n]|} \sum_{k\in \mathcal{U}\cap[1,n]} f(\tau^k x) =\int_G f d\mathbf{m} .
	\end{equation}
	Indeed, we have only to check \eqref{eq:tau-powers} for $f=\gamma \in \widehat{G}\setminus \{1\}$. As $\tau^*$ is monic, the functions $f\circ \tau^k = {\tau^*}^k \gamma$
	are distinct characters (the ergodicity is used here), hence orthogonal. Thus \eqref{eq:tau-powers} follows, because a law of large numbers applies. 
	
	Next, to establish existence of the limit 
    in \eqref{eq:tau-powers} for every $f\in L^1(G)$(regardless of its value),
by the Banach principle, we only need to establish a maximal inequality for the averages
	       $$
	      T_nf(x):=  \frac{1}{|\mathcal{U}\cap [1, n]|} \sum_{k\in \mathcal{U}\cap[1,n]} f(\tau^k x)
	       $$
	for all $f\in L^1(\mathbb{T})$. Since $\mathcal{U}$ has positive lower density, as pointed out by Conze \cite[Lemme 4]{Conze1973}, it is easy to deduce such an inequality from the 
	classical Hopf maximal inequality, just because assuming $f\ge 0$ we have
	$$
	    T_nf(x) \le   \frac{n}{|\mathcal{U}\cap [1, n]|} \cdot \frac{1}{n}\sum_{k\in [1,n]} f(\tau^k x) \le C \cdot \frac{1}{n}\sum_{k\in [1,n]} f(\tau^k x)
	$$
	where $C>0$ is a constant depending on the positive lower density.
	
	Now let us determine the limit in \eqref{eq:tau-powers}, which is denoted by $L_f(x)$ for the moment.  For any $\epsilon >0$, there exists a trigonometric polynomial $g_\epsilon$ such that $\|f-g_\epsilon\|_1<\epsilon$.
	By what we have just proved, we have
	\begin{equation}\label{eq:limit}
	     L_f(x) = L_{f-g_\epsilon}(x) + \int g_\epsilon d\mathbf{m} \ \ \ a.e.
	     \end{equation}
	     Also we have clearly
	     \begin{equation}\label{eq:limit2}
	      \left|\int g_\epsilon d\mathbf{m} -\int fd\mathbf{m}\right|<\epsilon.
	     \end{equation}
	 Now let us deal with $L_{f-g_\epsilon}$. By Fatou lemma, we get
	     $
	     \|L_{f-g_\epsilon}\|_1\le \|f-g_\epsilon\|_1.
	$
	Then,
	by the Markov inequality, we obtain the estimate
	$$
	         \mathbf{m} \{x:|L_{f-g_{\epsilon}}(x)|>\sqrt{\epsilon}\}\le \sqrt{\epsilon}.
	$$
	Then the Borel-Cantelli lemma implies  $|L_{f-g_{n^{-4}}}(x)|=O(n^{-2})$ as $n\to \infty$ for almost all $x$. Thus from \eqref{eq:limit} and \eqref{eq:limit2}, we get
	$$
	    a.s. \ \ \ L_f(x)= \int f d \mathbf{m}+ O(n^{-2}), 
	$$
	which gives us $L_f(x)= \int f d\mathbf{m}$ a.e.
	\end{proof}
	

	The question is more subtle when $\mathcal{U}$ has zero density. Clearly, the concept of Khintchin sequence is related to  the concept of universal good sequence in ergodic theory (see \cite{RW1995}). 
	The works on universal good sequences and universal bad sequences provide us information about $\tau^\mathcal{U}$.
	For example,  $\tau^{\{n^2\}}$ is a $L^p$-Khintchin sequence for $p>1$ by a result of Bourgain \cite{Bourgain1989} and is not a $L^1$-Khintchin sequence by a result of Buczolich-Mauldin \cite{BM2010},
	as is already mentioned in the special case $\tau(x) = 2 x \mod 1$.
    
	We can also consider a random set of powers instead of $\mathcal{U}$.
	Let $\mathbf{p}=(p_n)_{n\ge 1}$ with $0\le p_n<1$ for $n\ge 1$ such that $\sum_{n=1}^\infty p_n =\infty$. Let $(\xi_n)$ be a $\mathbf{p}$-Bernoulli sequence, i.e. a sequence of independent $0$-$1$ valued random variables 
	such that 
	$$
	    P(\xi_n=1)=p_n.
	$$
	We consider the random set $\mathcal{R}$ of integers, also denoted $\mathcal{R}_\mathbf{p}$ and called $\mathbf{p}$-Bernoulli set,   defined by
	$$
	      \mathcal{R} :=\mathcal{R}_\mathbf{p}:= \{n \ge 1:   \xi_n =1\}=\{n_1(\omega)< n_2(\omega)< \cdots\},
	$$
	and the set of random  powers of $\tau$:
	 $$
	 \tau^\mathcal{R} :=\tau^{\mathcal{R}_\mathbf{p}}:= \{\tau^{n_1(\omega)}, \tau^{n_2(\omega)},  \cdots\}.
	 $$
	 Bourgain \cite{Bourgain1988} proved that when $p_n = \frac{(\log \log n)^B}{n}$ ($n\ge 3$) with $B>1$. Then  
	 $ \tau^\mathcal{R}$ is almost surely a $L^p$-Khintchin set for $p> 1+ \frac{1}{B}$. 
     Remark that $ \tau^\mathcal{R}$ has zero density
     with respect to $\tau^{\mathbb{N}}$.
	 
	  We have the following corollary of Theorem \ref{thm:positivedensity}.
	 
	 \begin{cor} Almost surely $\tau^\mathcal{R} $ is a $L^1$-Khintchin sequence under the assumption
	 $$
	      \varliminf_{n\to \infty} \frac{p_1+\cdots +p_n}{n} >0.
	$$
	 \end{cor}
	 \begin{proof} We only need to show that $\mathcal{R}$ has almost surely a positive lower density. To this end, observe that
	 $$
	      \frac{|\mathcal{R}\cap [1, n]|}{n}=\frac{1}{n} \sum_{k=1}^n \xi_k =\frac{\sum_{k=1}^n p_k}{n} \times \frac{\sum_{k=1}^n \xi_k}{\sum_{k=1}^n p_k} .
	$$
	Then we get that 
	$$
	    a.s.  \quad \varliminf_{n\to \infty}  \frac{|\mathcal{R}\cap [1, n]|}{n}=\varliminf_{n\to \infty}  \frac{\sum_{k=1}^n p_k}{n} >0,
	$$
	because 
	\begin{equation}\label{eq:LLN}
	     a.s. \quad  \lim_{n\to \infty}\frac{\sum_{k=1}^n \xi_k}{\sum_{k=1}^n p_k} =1. 
	\end{equation}
	Indeed, let $s_n =p_1+p_2+\cdots +p_n$ and consider the martingale
	$$
	        M_n = \sum_{k=1}^n \frac{\xi_k -p_k}{s_k},
	$$
	which is $L^2$-bounded
	$$
	     \mathbb{E}(M_n^2) = \sum_{k=1}^n \frac{p_k(1-p_k)}{s_k^2} \le  \sum_{k=1}^\infty \frac{p_k}{s_k^2}<\infty. 
	$$
	Then, by Doob's convergence theorem, the series  $\sum_{k=1}^\infty \frac{\xi_k -p_k}{s_k}$ converges almost surely. Then \eqref{eq:LLN} follows, from the Kronecker lemma.
	 \end{proof}

     \subsection{ud-sequences of epimorphisms and products of epimorphisms}	

Any sequence of distinct epimorphisms on $\mathbb{T}$ is a ud-sequence. But it is not the case on  general groups. For example, as we have already mentioned, no sequence of 
epimorphisms on $\mathbb{Z}_p$, which is considered as an additive group, is a ud-sequence.

Here is a condition for a sequence of epimorphisms $\{\tau_n\}$ to be a ud-sequence.
	
	\begin{prop} \label{prop:ud}If $\{\tau_n\}\subset {\rm Epi}(G)$ is a sequence of epimorphisms such that 
	$\tau_n^*\gamma\not=\tau^*_m\gamma$ for any $\gamma\in \widehat{G}\setminus\{1\}$ and any integers
	$n\not=m$, then $\{\tau_n\}$ is a ud-sequence and any reordering of $\{\tau_n\}$ is still a ud-sequence. Consequently, for any continuous function 
    $f\in C(G)$, we have
    $$
    a.e. \quad \lim_{N\to \infty} \frac{1}{N}
    \sum_{n=1}^N f(\tau_n x) = \int_{G} f  d\mathbf{m} . 
    $$
	\end{prop}

	\begin{proof} By Weyl's criterion, it suffices to show that for any $\gamma\in \widehat{G}\setminus\{1\}$ we have
		$$
	       \lim_{N\to \infty} \frac{1}{N}\sum_{n=0}^{N-1} \gamma(\tau_n x)  = \lim_{N\to \infty} \frac{1}{N}\sum_{n=0}^{N-1} \tau_n^*\gamma(x) = 0 \ \ \ \ \mathbf{m}{\rm -}a.e.
	$$
	The assumption means that $\tau_n^*\gamma$'s for different $n$'s are distinct group characters. Then they are orthonormal functions in 
	$L^2(G)$ and satisfies a law of large numbers (for example, we can apply  the well known Davenport-Erd\"{o}s-LeVeque theorem \cite{DEL1963}). Having proved that $(\tau_n)$
    is a ud-sequence, we get immediately  equality by Weyl's criterion.
	\end{proof}

	Given a sequence of epimorphisms $\{A_n\}$ on $G$, we construct a new sequence of epimorphisms 
    by  composition
	$$
	    \tau_n =A_{n-1} \cdots A_1A_0.
	$$
	Here is a sufficient condition on $\{A_n\}$ such that  $\{A_{n-1}\cdots A_1A_0 x\}_{n\ge 1}$ is uniformly distributed for almost every $x\in G$, as a consequence of Proposition \ref{prop:ud}: 
	for any integers  $0\le n\le m<\infty$ and  any non trivial character $\gamma\in \widehat{G}\setminus\{1\}$, 
	$\gamma(A_{m-1}\cdots A_{n+1}x)$ and $\gamma(x)$ are different characters.	 
    \medskip
	
	But the above condition is not practical, because products $A_{m-1}\cdots A_{n+1}$ are involved.
	When $G=\mathbb{T}^d$,  we present a  more practical sufficient condition. We say that a $d\times d$ matrix $A$ is {\em expanding} if  
	\begin{equation}\label{eq:exp}
	 \quad \inf_{\|u\|=1} \|A u\|>1
	\end{equation}
	where $\|u\|$ denotes the Euclidean norm of $u$. Recall that the symmetric matrix $A^T\!A$ has nonnegative eigenvalues that can be
	ordered as $0\le \lambda_1\le \lambda_2\le \cdots \le \lambda_d$ and $\sqrt{\lambda_j}$ ($1\le j\le d$) are called the {\em singular values} of $A$. 
	
	\begin{lem} \label{lem:exp}Let $A$ be a real matrix. The following are equivalent
	\begin{itemize}
	\item[(1)]  $A$ is expanding;
	\item[(2)]  $A^T$ is expanding;
	\item[(3)] the singular values of $A$ are larger than $1$. 
	\end{itemize}
	\begin{proof} The equivalence of (1) and (3) follows from the equality
	$$
	       \|A u\|^2 = \langle Au, Au\rangle =  \langle A^T\!Au, u\rangle 
	$$
	and the fact that the eigenvectors of $A^T\!A$ are orthogonal. Indeed, assume $A^T\!A u_j =\lambda_j u_j$ with $\|u_j\|=1$.
	For an arbitrary unit vector $u$, which can be written as  $u =\sum_{j=1}^d \alpha_j u_j$ with $\sum_{j=1}^d \alpha_j^2=1$, we have
	$$
	 \|A u\|^2 = \sum_{j=1}^d \lambda_j \alpha_j^2. 
	$$
	
	The equivalence of (1) and (2) follows from the equivalence of (1) and (3), and the fact that
	$ A^T\!A$ and $AA^T$ have the same eigenvalues. 
	\end{proof}
	\end{lem}

	The condition in the following theorem is practical, because it is put on individual matrix
    and is easy to check.  
	
	\begin{thm} \label{prop:Exp-ud} Let  $\{A_n\}$  be a sequence of expanding endomorphisms on $\mathbb{T}^d$. Then  the sequence $\{A_{n-1}\cdots A_1A_0\}$ is a ud-sequence.
    Hence, for any 
    $f\in C(\mathbb{T}^d)$, we have
    $$
    a.e. \quad \lim_{N\to \infty} \frac{1}{N}
    \sum_{n=1}^N f(A_n \cdots A_2A_1 x) = \int_{\mathbb{T}^d} f d\mathbf{m} . 
    $$
	
	\end{thm}
	\begin{proof} Let $\tau_n =A_{n-1} \cdots A_1A_0$. It suffices to check the condition
	$\tau_n^*\gamma \not= \tau^*_m \gamma$ in Proposition \ref{prop:ud}., which can be restated as follows
	$$
	 \forall v \in \mathbb{Z}^d\setminus \{0\}, \forall n<m, \quad v A_{n-1} \cdots A_1A_0\not= v A_{m-1} \cdots A_1A_0.
	$$
	Since the matrices $A_k$'s are invertible, we have only to check $v\not=  v A_{m-1} \cdots A_{n+1}A_n$.
	That is true: first notice that, by Lemma \ref{lem:exp}, the transpose $A_n^T$ is expanding so
	$$
	     \lambda(A_n):= \inf_{\|v\|=1}\|v A_n\| >1.
	$$ 
	Hence we have
	$$
	  \|v A_{m-1} \cdots A_{n+1}A_n\| \ge \|v A_{m-1} \cdots A_{n+1}\| \lambda(A_n) \ge  \|v\|\lambda(A_{m-1})\cdots \lambda(A_n)>\|v\|.
	$$
	\end{proof}
	
	Theorem \ref{prop:Exp-ud} affirms that   the products of expansive endomorphisms on $\mathbb{T}^d$ are ud-sequences.  But
	Theorem \ref{prop:Exp-ud} can not apply to automorphisms on $\mathbb{T}^d$ ($d\ge 2$), which are not expanding. Some sequences of automorphisms on $\mathbb{T}^d$ with $d\ge 2$ may be
	ud-sequences and some are not, as the following examples show. First remark that there are only two automorphisms on $\mathbb{T}$: $x \mapsto x$ and $x \mapsto -x$. So, there is no sequence of distinct automorphisms  on $\mathbb{T}$, hence no question of ud-sequences of automorphisms.

	\medskip
	{\bf Example 1.} Consider the matrices $B_n\in SL_2(\mathbb{Z})$ defined by
	$$
	       B_n =\begin{pmatrix}
	           b_n & 1\\
	           1 & 0
	       \end{pmatrix} \quad {\rm with} \ \ b_n \in \mathbb{Z}.
	$$
    The characteristic polynomial of $B_n$ is $\lambda^2 -b_n \lambda -1$ and the eigenvalues are $\frac{1}{2}\big(b_n \pm \sqrt{b_n^2 +4}\big)$.
	Notice that $\det B_n =-1$. If $b_n=0$, $B_n$ has $\pm 1$ as eigenvalues. If $b_n\not=0$, one eigenvalue is larger than 1 in absolute value and the other less than 1, $B_n$ is then
	hyperbolic, hence  ergodic. Assume $b_n$'s are distinct,  then $B_n$ are distinct automorphisms of $\mathbb{T}^2$. 
    But the condition in Proposition \ref{prop:ud} is not satisfied. Indeed,
    we have 
    $$
    B_nB_m^{-1} = \begin{pmatrix}
	           b_n & 1\\
	           1 & 0
	       \end{pmatrix} \begin{pmatrix}
	           0 & 1\\
	           1 & -b_m
	       \end{pmatrix}
           =\begin{pmatrix}
	           1 & b_n-b_m\\
	           0 & 1
	       \end{pmatrix}.
    $$
    It follows  that  $B_nB_m^{-1}$'s has $(v_1, v_2)$ with $v_1=0$ as  left fixed point in $\mathbb{Z}^2\setminus\{(0,0)\}$, just because
	$$
	     (v_1, v_2) B_nB_m^{-1} = (v_1, (b_n-b_m) v_1 + v_2).
	$$
	 It is easy to check, by Weyl criterion,    that $\{B_n\}$ is not a ud-sequence. In fact, for any point 
     $\mathbf{x}=(x_1, x_2)^T$, i.e.
	$\{B_n \mathbf{x}\}$
	is not uniformly distributed.  Indeed, $B_n \mathbf{x} =  (b_n x_1+x_2, x_1)^T$. For the  character
    $\gamma(\mathbf{x}) = e^{2\pi i x_2}$, we have
    $$
    \frac{1}{N}\sum_{n=1}^N \gamma (B_n \mathbf{x}) = e^{2\pi i x_1}
    $$
    which does not tend to $0$.
    Actually, the orbit $\{B_n \mathbf{x}\}$ of $\mathbf{x}$ remains  in the one-dimensional torus $\mathbb{T}\times \{x_1\}$. 
	
	\medskip
    Now we present a ud-sequence of automorphisms.
    
	{\bf Example 2.} Now consider the matrices $B_n\in SL_2(\mathbb{Z})$ defined by
	$$
	       B_n =\begin{pmatrix}
	           b_n & b_n^2 -1\\
	           0 & b_n
	       \end{pmatrix} \quad {\rm with} \ \ b_n \in \mathbb{Z}\setminus\{0\}.
	$$
	Notice that $\det B_n =1$. 
	Assume $b_n$ are distinct.  If $n \not= m$, then $B_n$ and $ B_m$ have no left fixed point in $\mathbb{Z}^2\setminus\{(0,0)\}$, just because
	$$
	     (v_1, v_2) B_n = (v_1 b_n, v_1 (b_n^2-1)+ b_n v_2).
	$$
	So, if $b_n$'s are distinct,  $\{B_n\}$ is a ud-sequence, i.e.
	$\{B_n x\}$
	is uniformly distributed for almost every $x \in \mathbb{T}^2$.  
    
    Let us consider the special cases (corresponding to $b_n=n$): 
	$$
	   B_n' =\begin{pmatrix}
	           n & n^2 -1\\
	           0 & n
	       \end{pmatrix}. 
	$$
	We claim  that $\{B_n'\}$ is not a $L^\infty$-Khintchin sequence. Indeed, if $g$ is a counter-example of Marstrand ensuring that $\{n\} \subset \mathbb{Z}$ is not a
	$L^\infty$-Khintchin sequence, and if we consider $f(x_1,x_2)=g(x_2) \in L^\infty(\mathbb{T}^2)$, we get 
	$$
	    \frac{1}{N} \sum_{n=1}^N f(B_n' x) = \frac{1}{N} \sum_{n=1}^N g(nx_2),
	$$  
	which doesn't tend  almost everywhere to $\int g(x) dx =\int f(x_1,x_2) dx_1 dx_2$. 
    \medskip

	To construct a ud-sequence or better a Khintchin sequence $\{\tau_n\}$, as we did for Khintchin sequence of integers, we  first select  a sequence of epimorphisms $\{A_n\}$ and then 
	consider the products $\tau_n =A_{n-1}\cdots A_1 A_0$. 
    
    Here is an example.  Let $\{t_n\} \in \{0, 1\}^\infty$ be the Thue-Morse sequence
	and let $C_0, C_1$ be two given epimorphisms. We select $A_n$ from $\{C_0, C_1\}$ by $\{t_n\} $:
	         $
	              A_n =C_{t_n}.
	         $    

             \begin{prop}
             Let  $\{t_n\} \in \{0, 1\}^\infty$ be the Thue-Morse sequence.
                 If two epimorphisms $C_0$ and $C_1$ on a compact abelian group commute each other  and  if the product endomorphism $C_0C_1$ is ergodic, then the sequence
             $(C_{t_1}C_{t_2}\cdots C_{t_n})$ is a Khintchin sequence.
             \end{prop}
	          The proof of this result is the same as that of Theorem \ref{thm:TM}.

             We don't know what happens when $C_0$ and $C_1$ are not commutative. In particular, we have the  following open question. 
             If we assume that $C_0$ and $C_1$ are expansive matrices but not commutative,  we know that $C_{t_n}\cdots C_{t_2}C_{t_1}$
	         is a ud-sequence (see Proposition \ref{prop:Exp-ud}). 

             \begin{Question}
             Assume that $C_0$ and $C_1$ are expansive but not commutative matrices. 
             Is the multiplicative Thue-Morse sequence $C_{t_n}\cdots C_{t_2}C_{t_1}$  a Khintchin sequence ?
             We also ask the question for other primitive sequences.
             \end{Question}

	\section{Skew products of endomorphisms on groups}\label{sect:skewproduct}

In this section,  we present a random construction,
using the skew products. In other words, we randomize our multiplicative procedure in order to produce Khintchin sequences. We find a condition, called Fourier tightness, ensuring that the skew product dynamics shares the same properties as the basis dynamical systems. These properties include ergodicity, weakly mixing and strongly mixing. 

\subsection{Skew products}

	The skew product of measure preserving transformations was introduced in its generality by Anzai \cite{Anzai1951} in 1951. Let us first recall the definition. 
	Let $(\Omega, \mathcal{A}, \mu, \sigma)$ be a measure-preserving dynamical system. Let $(X, \mathcal{B}, \nu)$ be another probability space. On the other hand, to 
	each $\omega \in \Omega$  is associated a measurable transformation $\psi_\omega: X\to X$ which preserves the probability measure $\nu$. 
	Suppose that $(\omega, x) \mapsto \psi_\omega (x)$ is a measurable mapping from $\Omega \times X$ into $X$. The {\em skew product transformation} of $\{\psi_\omega\}_{\omega\in \Omega}$
	based on $(\Omega, \mathcal{A}, \mu, \sigma)$,
	denoted $S: \Omega \times X \to \Omega\times X$, is defined by
	$$
	       S(\omega, x) = (\sigma \omega, \psi_\omega(x)).
	$$
		Then $S$ preserves the product measure   $\mu \otimes \nu$ 
		of $\mu$ and $\nu$.
		The measure preserving dynamical system $(\Omega \times X, \mathcal{A}\times \mathcal{B}, \mu \otimes \nu, S)$ is called the {\em skew product} of measure preserving 
		systems $(X, \mathcal{B}, \nu, \psi_\omega)$ based on $(\Omega, \mathcal{A}, \mu, \sigma)$. If $\psi_\omega$ is constantly equal to a fixed transformation $\psi$,
		we get the direct product of $(\Omega, \mathcal{A}, \mu, \sigma)$  and $(X, \mathcal{B}, \nu,\psi)$.
		
		In this paper, we specialize to the  case where $X=G$ is a compact abelian group,  $\nu=\mathbf{m}$ (sometimes denoted $dx$) is the Haar probability measure on the Borel field $\mathcal{B}$ of $G$ and $\psi_\omega$'s
		are endomorphisms on $G$:
		 $$
		     \psi_\omega = A_\omega x
		 $$
		 where  $A_\omega$ is a surjective endomorphism of $G$ (we can also consider compact non abelian groups, but they are not concerned in the present paper).


		Our model can be defined as follows.  
	Let $\Omega ={\rm Epi}(G)^\mathbb{N}$  and $\sigma$ be the shift on $\Omega$ defined by $(\omega_n)_{n\ge 0}\mapsto (\omega_{n+1})_{n\ge 0}$.
	Then the skew product transformation  $S: \Omega \times G \to \Omega \times G $ is defined by 
	$$
	   S(\omega, x) =(\sigma \omega, \omega_0 x).
	$$
	Recall that if $\mu$ is a $\sigma$-invariant measure, then $\mu\times \mathbf{m}$ is $S$-invariant.

	We will study the ergodicity, weakly mixing property and strong mixing property of our skew products of epimorphisms.
	For simplicity, we will say that a probability measure $\nu$ is $T$-ergodic if $\nu$ is invariant under the map $T$ and is ergodic. Similarly, we talk about
	weakly $T$-mixing and strong $T$-mixing.

	\subsection{Fourier tightness}
	
	We first discuss an assumption that we will make.
	As usual, for $f\in L^1(\mu\times \mathbf{m} )$, we denote $f\circ S$ by $Sf$. Thus for a character $\gamma \in \widehat{G}$, we have $S^n \gamma (\omega, x)= \gamma(\omega_{n-1}\cdots \omega_1 \omega_0 x)$, which could be considered as a random character of $G$. It is nothing but $(\omega_{n-1} \cdots \omega_1\omega_0)^*(\gamma)$. We will  make use of the following
	 asssumption: almost surely,  for any finite set of characters $\Gamma \subset \widehat{G}$, there exists an integer $m_0=m_0(\omega,\Gamma)$ such that
	$$
	     \forall m\ge m_0, \quad    (\omega_{m-1}\omega_{m-2}\cdots \omega_1\omega_0)^*(\widehat{G}\setminus\{1\}) \subset  \Gamma^c.
	        $$
	        That is to say, for large $m$, the adjoint $(\omega_{m-1}\omega_{m-2}\cdots \omega_1\omega_0)^*$ pushes $\widehat{G}\setminus\{1\}$ outside $\Gamma$.
	  If this assumption is verified, 
	        we say that the skew product of epimorphisms is {\em Fourier tight}.
	        
	        On the group $\mathbb{T}$,  ${\rm Epi}(\mathbb{T})$ is identified with $\mathbb{Z}^*:=\mathbb{Z}\setminus\{0\}$ so that $\Omega = (\mathbb{Z}^*)^{\mathbb{N}}$. 
            For an integer $a$, denote by $[a]$ the set of all sequences in $\Omega$ whose first term is $a$.
	        If $\mu$ is a $\sigma$-invariant ergodic measure such that $\mu([a])>0$ for some integer $a$ with $|a|\ge 2$, then our  skew product of epimorphisms is Fourier tight.
	        Indeed, if we denote \[N_n(\omega, a)=\#\{0\le k\le n-1: \omega_k =a\},\] Birkhoff's ergodic theorem implies that
	        $$
	            a.s. \quad \lim_{n\to \infty}\frac{N_n(\omega, a)}{n} = \lim_{n\to \infty}\frac{1}{n} \sum_{k=0}^{n-1} 1_{[a]}(\sigma^k \omega)=\mu([a])>0.
	        $$
	        So, almost surely, there exists an integer $n(\omega)$ such that $N_n(\omega, a)> \frac{\mu([a])}{2} n$ for $n \ge n(\omega)$. It follows that
	        \begin{equation}\label{eq:FT1}
	             |\omega_{n-1}\cdots \omega_1\omega_0|\ge |a|^{\frac{\mu([a]) n}{2}}> 2^{\frac{\mu([a]) n}{2}}
	        \end{equation}
	        which implies  that $(\omega_{n-1}\cdots \omega_1\omega_0)^*(\mathbb{Z}^*)$  is outside the interval $[-2^{\frac{\mu([a]) n}{2}}, 2^{\frac{\mu([a]) n}{2}}]$.
            \medskip
	        
	        Notice that there are only two automorphisms corresponding to $-1$ and $1$ on the group $\mathbb{T}$. All other epimorphisms are "expanding". The above 
	        proof of Fourier tightness benefits from this expanding property.  On the group $\mathbb{T}^d$ with $d\ge 2$, there are 
	        many epimorphisms which are not expanding. 
	        The same argument as above show that  skew product of epimorphisms on $\mathbb{T}^d$ is Fourier tight if the $\sigma$-invariant measure $\mu$
	        is ergodic and positively charges a set of   expanding matrices. 

	    \subsection{Ergodicity}     
	\begin{thm}\label{thm:ergodic} Let $G$ be a compact abelian group with countable dual group $\widehat{G}$.
	Suppose that the skew product of epimorphisms is Fourier tight. Then
	 $\mu\otimes \mathbf{m}$ is $S$-ergodic iff
	$\mu$ is  $\sigma$-ergodic.
	\end{thm}
	
	\begin{proof} Since $\sigma$ is a factor of $S$, the $S$-ergodicity implies the $\sigma$-ergodicity. We should only to prove the converse proposition.
	
	Assume that $\mu$ is $\sigma$-ergodic. 
	   Let $f \in L^2(\mu \times \mathbf{m})$ be a $S$-invariant function. That is to say, 
	   \begin{equation}\label{eq:inv}
	         f(\sigma \omega, \omega_0 x) = f(\omega, x)   \quad \mu\times \mathbf{m}\!-\!a.e.
	   \end{equation}
	   We are going to show that $f$ is $\mu\times \mathbf{m}$-almost everywhere constant. 
	   
	   Denote $f_\omega(x) =  f(\omega, x)$. By Fubini theorem, for $\mu$-almost all $\omega$, $f_\omega$ is a well defined function in $L^2(\mathbf{m})$. First notice that 
	   since the surjective endormorphism $x\mapsto \omega_0 x$ preserves $\mathbf{m}$, the invariance (\ref{eq:inv}) implies 
	   $$
	        \|f_{\sigma \omega }(\cdot)\|_{L^2(\mathbf{m})}=  \|f_{\omega}(\cdot)\|_{L^2(\mathbf{m})}.
	   $$
	   Thus the norm function $\omega \mapsto  \|f_{\omega}(\cdot)\|_{L^2(\mathbf{m})}$ is  $\sigma$-invariant, then   $\mu$-a.e. constant by the ergodicity of $\mu$.
	   By the same reason, the integral function $\omega\mapsto \int f_\omega d\mathbf{m}$ is also $\sigma$-invariant and then constant.  So, without loss of generality, we can assume that 
	   $\int f_\omega d\mathbf{m} =0$. 
	   
	   Now consider the Fourier series of $f_\omega(\cdot)$
	   $$
	          f_\omega(x) =\sum_{\gamma \in \widehat{G}\setminus\{1\}} \widehat{f}_\omega(\gamma) \gamma(x).  
	   $$ 
	   Then the invariance (\ref{eq:inv}) is translated into 
	   \begin{equation}\label{eq:inv2}
	         \sum_{\gamma \in \widehat{G}\setminus\{1\}} \widehat{f}_{\sigma\omega}(\gamma) \omega_0^* (\gamma) = \sum_{\gamma \in \widehat{G}\setminus\{1\}} \widehat{f}_\omega(\gamma) \gamma.
	   \end{equation}
	   In other words, we have
	    \begin{equation}\label{eq:inv3}
	     \forall \gamma\in \widehat{G}\setminus\{1\}, \ \ \   \widehat{f}_{\sigma\omega}(\gamma) = \widehat{f}_\omega(\omega_0^*(\gamma)).
	     \end{equation}
	     By iteration,  for all $\gamma\in \widehat{G}\setminus\{1\}$ and all $n\ge 1$ we get 
	     \begin{equation}\label{eq:FourierRec}
	         \widehat{f}_{\sigma^n\omega}(\gamma) =\widehat{f}_\omega(\omega_{0}^*\omega_1^*\cdots \omega_{n-1}^*(\gamma))  = \widehat{f}_\omega((\omega_{n-1}\cdots \omega_1\omega_0)^*(\gamma)).
	     \end{equation}
	     So, by the invariance of the norm $\|f_{\omega}\|_{L^2(\mathbf{m})}$ and the  recursive formula \eqref{eq:FourierRec},  we get
	     $$
	      \forall n\ge 1, \quad  \|f_{\omega}\|_{L^2(\mathbf{m})}^2 = \|f_{\sigma^n \omega}\|_{L^2(\mathbf{m})}^2
	       = \sum_{\gamma \in \widehat{G}\setminus\{1\}} |\widehat{f}_{\omega}((\omega_{n-1}\cdots \omega_1\omega_0)^*(\gamma))|^2.
	     $$
	     Let $\Gamma$ be any finite subset of $\widehat{G}$.
	     The assumption of Fourier tightness implies  that almost surely there exists $n(\omega)$ such that
	     $$
	       \forall n\ge n(\omega), \ \ \  \|f_{\omega}\|_{L^2(\mathbf{m})}^2  \le \sum_{\gamma \in \Gamma^c}|\widehat{f}_\omega(\gamma)|^2. 
	        	     $$
		This,  together with Parseval's identity, implies that $\|f_{\omega}\|_{L^2(\mathbf{m})}^2=0$ a.s., so $f=0$.
	\end{proof}

	\subsection{Weakly mixing property and strong mixing property}

	If $\mu\times \mathbf{m}$ is $S$-weakly mixing, $\mu$ must be $\sigma$-weakly mixing. 
	We have the following  converse.
	
	\begin{thm} \label{thm:w-mixing} Let $G$ be a compact abelian group with countable dual group $\widehat{G}$.
	Suppose that the skew product of epimorphisms is Fourier tight.
	If $\mu$ is weakly $\sigma$-mixing, then $\mu\times \mathbf{m}$ is weakly $S$-mixing.
	\end{thm}
	
	\begin{proof} 
	That $\mu\times \mathbf{m}$ is $S$-weakly mixing means that $(\mu\times \mathbf{m}) \times (\mu\times \mathbf{m})$ is $S\times S$-ergodic (\cite{Walters1982}). 
	Recall that $S\times S: (\Omega\times G) \times  (\Omega\times G)$ is defined by 
	$$
	    (\omega', x;\  \omega'', y) \mapsto (\sigma\omega', \omega'_0 x;\ \sigma\omega'', \omega''_0 y).
	$$
	The product $S\times S$ can be viewed as a skew product of endomorphisms of $G\times G$ over the shift $\sigma\times \sigma$ on $\Omega\times \Omega$.
	In other words, $S\times S$ can be defined on $(\Omega \times \Omega) \times (G\times G)$ by
	$$
	    (S\times S)(\omega', \omega''; x,y) \mapsto (\sigma\omega', \sigma\omega''; \omega_0'x, \omega''_0 y).
	$$
	Since $\mu$ is assumed weakly $\sigma$-mixing, $\mu\times \mu$ is $\sigma\times \sigma$-ergodic. By Theorem \ref{thm:ergodic}, 
	$(\mu\times \mathbf{m}) \times (\mu\times \mathbf{m})$ is $S\times S$-ergodic.   	
	\end{proof}
	
	Theorem \ref{thm:w-mixing} can be also deduced from Theorem \ref{thm:eigenvalue} which will be proved  in the next subsection. 
	
	The strong mixing property is also preserved when we  extend $\sigma$ to $S$.
	
	\begin{thm} \label{thm:mixing} Let $G$ be a compact abelian group with countable dual group $\widehat{G}$.
	Suppose that the skew product of epimorphisms is Fourier tight. If $\mu$ is $\sigma$-mixing, then $\mu\otimes \mathbf{m}$ is $S$-mixing.
	\end{thm}
	
	\begin{proof} Recall that $\mu\otimes \mathbf{m}$ is $S$-mixing means that for any $F, G\in L^2(\mu \otimes \mathbf{m})$, we have
	$$\lim_{n\to \infty }\int F\cdot G\circ S^n d\mu \otimes \mathbf{m}= \int F d\mu \otimes \mathbf{m} \cdot \int G d\mu \otimes \mathbf{m}.$$
	It suffices to justify this for $F$ and $G$ belonging to a dense family of functions. So we can assume that
	$F(\omega,x) = f_1(\omega)f_2(x)$ and    $G(\omega,x) = g_1(\omega)g_2(x)$ with $f_1, g_1 \in L^\infty(\mu)$
	and $f_2, g_2\in L^\infty(\mathbf{m})$. We can even suppose that $f_2$ and $g_2$ are trigonometric polynomials.
	
	Let us first observe that
	\begin{equation}\label{eq:mix1}
	   \int F  \cdot  G \circ S^n  d\mu \otimes \mathbf{m} = \int_{\Omega} f_1(\omega) g_1(\sigma^n \omega) 
	   \left(
	   \int_{\mathbb{T}} f_2(x)g_2(\omega_{n-1}\cdots \omega_1\omega_0 x) d\mathbf{m}(x)
	   \right) d\mu(\omega).
	\end{equation}
	Since $\sigma$ is mixing, we have 
	\begin{equation}\label{eq:mixing}
	   \lim_{n\to \infty} \int_{\Omega} f_1(\omega) g_1(\sigma^n \omega) d\mu(\omega)= \int_{\Omega} f_1(\omega)d\mu(\omega) \int_\Omega g_1(\omega) d\mu(\omega).
	\end{equation}
	So, we are led to prove
	\begin{equation}\label{eq:F2}
	   \lim_{n\to \infty}   \int_{G} f_2(x)g_2(\omega_{n-1}\cdots \omega_1\omega_0 x) d\mathbf{m}(x) = \widehat{f}_2(1) \widehat{g}_2(1).
	\end{equation}
	
	If we develop $f_2$ by its Fourier series (a finite sum because we assume that $f_2$ is a polynomial), we get 
	\begin{eqnarray*}
	 \int_{G} f_2(x)g_2(\omega_{n-1}\cdots \omega_1\omega_0 x) d\mathbf{m}(x)
	 & = & \sum_{\gamma\in\widehat{G}} \widehat{f}_2(\gamma) \int \gamma(x) g_2(\omega_{n-1}\cdots \omega_1\omega_0 x)d\mathbf{m}(x).
	 \end{eqnarray*}
	 Now develop $g_2$ into Fourier series so that 
	 $$
	 g_2(\omega_{n-1}\cdots \omega_1\omega_0 x) = \sum_{\gamma' \in \widehat{G}}  \widehat{g}_2(\gamma')(\omega_{n-1}\cdots \omega_1\omega_0)^*(\gamma')(x).
	 $$
	 Putting this into the last equality and then integrating, we get
	 \begin{align}\label{eq:4.9}
	 \int_{G} f_2(x)g_2(\omega_{n-1}\cdots \omega_1\omega_0 x) d\mathbf{m}(x)
	 = \widehat{f}_2(1) \widehat{g}_2(1) + \sum_{\gamma' \in \widehat{G}\setminus \{1\}} \widehat{f}_2( (\omega_{n-1}\cdots \omega_1\omega_0)^* (\overline{\gamma'})) \widehat{g}_2(\gamma'). 
	 \end{align}
	Let $R_n(\omega)$ be the last sum.
	 By the Cauchy-Schwarz inequality, we get
	 $$
	    |R_n(\omega)| \le \|g_2\|_2 \sqrt{\sum_{\gamma' \in \widehat{G}\setminus \{1\}} |\widehat{f}_2( (\omega_{n-1}\cdots \omega_1\omega_0)^* (\overline{\gamma'}))|^2}
	 $$
	 which  tends to zero almost surely by the Fourier tightness. Thus \eqref{eq:F2} is proved. 
	 
	 Using  \eqref{eq:4.9}, write \eqref{eq:mix1} as follows
	 $$
	      \int F  \cdot  G \circ S^n  d\mu \otimes \mathbf{m} = \widehat{f}_2(1) \widehat{g}_2(1) \int_{\Omega} f_1(\omega) g_1(\sigma^n \omega) d\mu(\omega)
	      +  \int_{\Omega} f_1(\omega) g_1(\sigma^n \omega) R_n(\omega)d\mu(\omega).
	 $$
	 Observe that $R_n(\omega)$ bounded by $\|f_2\|_2 \|g_2\|_2$ for all $\omega$ so that 
	$ f_1(\omega) g_1(\sigma^n \omega) R_n(\omega)$ are a.s. bounded by the constant $\|f_1\|_\infty\|g_1\|_\infty\|f_2\|_2 \|g_2\|_2$. 
	Thus applying the Lebesgue dominated convergence theorem, we get 
	$$
	    \lim_{n\to \infty}\int_{\Omega} f_1(\omega) g_1(\sigma^n \omega) R_n(\omega)d\mu(\omega) =0.
	$$
	 So, we can really take limit in (\ref{eq:mix1}) by using \eqref{eq:mixing} in order to get
	 
\begin{align*}
	    \lim_{n\to \infty}   \int F d\mu \otimes \mathbf{m} \cdot \int G \circ S^n  d\mu \otimes \mathbf{m} =& \widehat{f}_2(1) \widehat{g}_2(1) 
	  \int f_1d\mu \int g_1d\mu \\
 =& \int F d\mu  \otimes \mathbf{m} \int G d\mu  \otimes \mathbf{m}.    
\end{align*}   
	\end{proof}

	\subsection{Eigenvalues}
	
	We are now going to prove   that the skew product dynamical system has the same eigenvalues as the basic dynamical system.
	In fact, the proof of Theorem \ref{thm:ergodic} can be almost repeated to prove the following result, which is stranger than Theorem \ref{thm:w-mixing} above
	concerning the weakly mixing property because a system is weakly mixing if and only if it has $1$ as simple eigenvalue and has no other eigenvalues.. 
	
	\begin{thm}\label{thm:eigenvalue} Suppose that $\mu$ is $\sigma$-ergodic. Then 
	the operator $F \mapsto F\circ S $ acting on $L^2(\mu\otimes\mathbf{m})$ has exactly the same eigenvalues as the operator
	$f \mapsto f\circ \sigma$  acting on $L^2(\mu)$.
	More precisely, $F\circ S = \lambda F$ with $F\in L^2(\mu \otimes \mathbf{m})$ iff $F(\omega, x) = f(\omega)$ with $f \in L^2)(\mathbf{m})$ such that $f \circ \sigma =\lambda f$.
	\end{thm}
	
	\begin{proof} It suffices to prove the ``only if" part, because the ``if" part is obvious. Suppose that the  eigen-equation  
	 \begin{equation}\label{eq:EV}
	 F(\sigma \omega, \omega_0 x)  = \lambda F(\omega, x)
	 \end{equation}
	 has a non-trivial solution   $F\in L^2(\mu \otimes \mathbf{m})$ for some $\lambda$ with $|\lambda|=1$. Since $\mu$ is assumed  ergodic, 
	 $\mu\otimes \mathbf{m}$ is $S$-ergodic (Theorem \ref{thm:ergodic}).  Denote $f_\omega(x) = F(\omega, x)$. From the eigen equation, we first get  the facts 
	 \begin{equation}\label{eq:EV1}
	 \|f_{\sigma \omega}\|_{L^2(\mathbf{m})} = \|f_\omega\|_{L^2(\mathbf{m})}; 
	 \quad |\widehat{f}_{\sigma^n\omega}(\gamma)| = |\widehat{f}_\omega((\omega_{n-1}\cdots \omega_1\omega_0)^* \gamma)|.
	 \end{equation}
	 Here we have used the fact $|\lambda|=1$ and the fact that the epimorphisms preserve the Haar measure. The second equality is similar to \eqref{eq:FourierRec})  
	 and it compares the absolute values of Fourier coefficients, not directly the Fourier coefficients. But it is sufficient for us.
	  As we did in the proof of Theorem \ref{thm:ergodic}, by (\ref{eq:EV1}) and the Parseval identity we deduce that
	  $$
	      \ \|f_\omega\|_{L^2(\mathbf{m})}^2 = \|f_{\sigma^n \omega}\|_{L^2(\mathbf{m})}^2
	       =  \left|\int f_\omega d\mathbf{m}\right|^2 + \sum_{\gamma \in \widehat{G}\setminus\{1\}} |\widehat{f}_{\omega}((\omega_{n-1}\cdots \omega_1\omega_0)^* \gamma )|^2
	     $$
	hold for all $m\ge 1$. 
	But the last sum tends to zero as $n\to \infty$, because of the Fourier tightness. Thus we get  
	$$
	        \int |f_{\omega}|^2 d\mathbf{m} =\left|\int f_\omega d\mathbf{m}\right|^2.
	$$
	This  implies that $f_\omega(\cdot)$ is constant. So, $F(\omega, x) =f(\omega)$ for some function $f$ defined on $\Omega$, which must be an eigenfunction of $\sigma$ associated to the eigenvalue $\lambda$
	(cf. (\ref{eq:EV})).
	\end{proof}

	\subsection{Random Khintchin sequences in weak sense}

	Recall that  $\Omega ={\rm Epi}(G)^\mathbb{N}$, $\sigma$ is the shift on $\Omega$ defined by $(\omega_n)_{n\ge 0}\mapsto (\omega_{n+1})_{n\ge 0}$
	and $\mu$ is a $\sigma$-invariant measure.
	The skew product of epimorphisms  $S: \Omega \times G \to \Omega \times G $ was defined by 
	$$
	   S(\omega, x) =(\sigma \omega, \omega_0 x).
	$$

	\begin{thm} \label{thm:WKS}Let $G$ be a compact abelian group with  dual group $\widehat{G}$.
	Suppose that the skew product of epimorphisms  is Fourier tight and that the measure 
	$\mu$ is  $\sigma$-ergodic.  Then for any $f\in L^1(G)$, $\mu$-almost surely
	\begin{equation}\label{eq:WKS}
	   \mathbf{m}{\rm -}a.e. \quad \lim_{n\to \infty}\frac{1}{n}\sum_{k=1}^n f(\omega_{k-1}\cdots \omega_1\omega_0 x) =\int_{G}  f d\mathbf{m}.
	\end{equation}
	\end{thm}

	\begin{proof} Let $F(\omega, x)=f(x)$.  It is a simple application of Theorem \ref{thm:ergodic} to the function $F\in L^1(\mu\times \mathbf{m})$, because
	$S^k(\omega, x)=(\sigma^k \omega, \omega_{k-1}\cdots \omega_1\omega_0 x)$. 
	\end{proof}
	
	\begin{Question}\label{qest: 8}
    Is the sequence $\{\omega_{k-1}\cdots \omega_1\omega_0\}$
    in Theorem \ref{thm:WKS} almost surely a Khintchin sequence on $G$\, ?
    \end{Question}

    Lacey, Petersen, Wierdl and Rudolph \cite{LPWR1994}
    studied random ergodic theorems with universally
representative sequences. This is related to Question \ref{qest: 8}. Unfortunately, the results in \cite{LPWR1994}
don't supply an answer to Question \ref{qest: 8}.
Actually transformations studied in \cite{LPWR1994} are invertible and our epimorphisms are not necessarily invertible.


	\end{document}